\documentclass[10pt,twoside]{article}

\usepackage{amssymb,amsmath,makeidx}
\usepackage[T1]{fontenc}
\usepackage{latexsym}
\usepackage{exscale}
\usepackage{pb-diagram}
\usepackage[dvips]{graphicx}
\usepackage{xr}
\usepackage{xfrac,faktor}
\usepackage{ulem}
\usepackage{xcolor,pgf,tikz}

\newcommand{\al}{\alpha}

\newcommand{\alpr}{{\alpha^\prime}}

\newcommand{\alnod}{{\al_0}}

\newcommand{\almmin}{{\al_{m-1}}}

\newcommand{\alm}{{\al_m}}

\newcommand{\alplus}{{\al^+}}

\newcommand{\be}{\beta}

\newcommand{\ga}{\gamma}

\newcommand{\gapr}{{\ga^\prime}}

\newcommand{\Ga}{\Gamma}
\newcommand{\Gapr}{{\Gamma^\prime}}
\newcommand{\de}{\delta}

\newcommand{\De}{\Delta}
\newcommand{\Depr}{{\Delta^\prime}}
\newcommand{\Deprnod}{{\Delta^\prime_0}}

\newcommand{\eps}{\varepsilon}

\newcommand{\epsn}{\varepsilon_0}

\newcommand{\la}{\lambda}

\newcommand{\om}{\omega}

\newcommand{\Om}{\Omega}

\newcommand{\Omi}{{\Omega_i}}
\newcommand{\Omipr}{{\Omega_i^\prime}}
\newcommand{\Omj}{{\Omega_j}}
\newcommand{\Omjpr}{{\Omega_j^\prime}}
\newcommand{\Omk}{{\Omega_k}}
\newcommand{\Oml}{{\Omega_l}}

\newcommand{\Omie}{\Om_{i+1}}
\newcommand{\Omiz}{\Om_{i+2}}
\newcommand{\Omje}{\Om_{j+1}}
\newcommand{\Omjz}{\Om_{j+2}}

\newcommand{\rhopr}{{\rho^\prime}}

\newcommand{\si}{\sigma}
\newcommand{\tht}{\vartheta}

\newcommand{\xipr}{{\xi^\prime}}

\newcommand{\ze}{\zeta}

\newcommand{\etapr}{{\eta^\prime}}

\renewcommand{\phi}{\varphi}

\newcommand{\N}{{\mathbb N}}

\newcommand{\Hz}{{\mathbb P}}

\newcommand{\Ez}{{\mathbb E}}

\newcommand{\CNF}{{\mathrm{\scriptscriptstyle{CNF}}}}
\newcommand{\ANF}{{\mathrm{\scriptscriptstyle{ANF}}}}

\newcommand{\NF}{{\mathrm{\scriptscriptstyle{NF}}}}

\newcommand{\Lim}{\mathrm{Lim}}

\newcommand{\sumend}{{\mathrm{end}}}

\newcommand{\thtn}{\tht_n}
\newcommand{\thtnod}{\tht_0}
\newcommand{\thte}{\tht_1}

\newcommand{\thti}{\tht_i}
\newcommand{\thtie}{\tht_{i+1}}

\newcommand{\thtj}{\tht_j}
\newcommand{\thtje}{\tht_{j+1}}

\newcommand{\T}{{\operatorname{T}}}

\newcommand{\Tt}{{\operatorname{T}^\tau}}

\newcommand{\stark}{{\ast_k}}

\newcommand{\mc}{{\operatorname{mc}}}

\newcommand{\Par}{\operatorname{Par}}

\newcommand{\PA}{{\operatorname{PA}}}

\newcommand{\pioneonecanod}{\Pi^1_1{\operatorname{-CA}_0}}
\newcommand{\pioneonetrnod}{\Pi^1_1{\operatorname{-TR}_0}}

\newcommand{\IDn}{{\operatorname{ID}_n}}
\newcommand{\IDnod}{{\operatorname{ID}_0}}

\newcommand{\qed}{\mbox{ }\hfill $\Box$\vspace{2ex}}

\newlength{\hilflh}

\newtheorem{theo}{Theorem}[section]
\newtheorem{cor}[theo]{Corollary}
\newtheorem{lem}[theo]{Lemma}
\newtheorem{defi}[theo]{Definition}
\newtheorem{prop}[theo]{Proposition}

\newtheorem{rmk}[theo]{Remark}

\newcommand{\lpr}{{l^\prime}}
\newcommand{\hpr}{{h^\prime}}

\newcommand{\Romannumeral}[1]{\uppercase\expandafter{\romannumeral #1\relax}}

\newcommand{\Si}{\Sigma}

\newcommand{\stari}{{\star_{i}}}

\newcommand{\starj}{{\star_{j}}}
\newcommand{\starje}{{\star_{j+1}}}

\newcommand{\chiome}{\chi^{\Om_1}}

\newcommand{\chiomie}{\chi^{\Omie}}

\newcommand{\Tomie}{\T^{\Omie}}
\newcommand{\Tomnod}{\T^{\Om_0}}
\newcommand{\Tomi}{\T^{\Omi}}

\newcommand{\chiomje}{\chi^{\Omje}}
\newcommand{\chiomjz}{\chi^{\Omjz}}

\newcommand{\domf}{\mathrm{d}}

\newcommand{\ual}{\underline{\al}}
\newcommand{\uual}{\uuline{\al}}

\newcommand{\ualpr}{{\ual^\prime}}

\newcommand{\ube}{\underline{\be}}

\newcommand{\uga}{\underline{\ga}}
\newcommand{\uGa}{\underline{\Ga}}

\newcommand{\Tcirc}{{\mathring{\T}}}

\newcommand{\Esc}{{\mathcal{E}}}
\newcommand{\Equok}{{\sfrac{\Esc}{k}}}
\newcommand{\Equol}{{\sfrac{\Esc}{l}}}
\newcommand{\Equokbig}{{\faktor{\Esc}{k}}}
\newcommand{\TcapOmquok}{{\sfrac{{\T\cap\Om_1}}{{\Ck}}}}
\newcommand{\TcapOmquoke}{{\sfrac{{\T\cap\Om_1}}{{\Cke}}}}
\newcommand{\Tquok}{{\sfrac{\T}{k}}}
\newcommand{\Tcircquok}{{\sfrac{\Tcirc}{k}}}

\newcommand{\Tscircquok}{{\sfrac{\Tcirc^s}{k}}}
\newcommand{\Tfcircquok}{{\sfrac{\Tcirc^f}{k}}}
\newcommand{\Tlcircquok}{{\sfrac{\Tcirc^l}{k}}}
\newcommand{\Tcircquokcls}{{\overline{{\sfrac{\Tcirc}{k}}}}}

\newcommand{\Tquokbig}{{\faktor{\T}{k}}}
\newcommand{\Tcircquokbig}{{\faktor{\Tcirc}{k}}}
\newcommand{\Tk}{{\Tcirc{\mbox{\raisebox{0.18ex}{$\scriptstyle [k]$}}}}}
\newcommand{\Tke}{{\Tcirc{\mbox{\raisebox{0.18ex}{$\scriptstyle [k+1]$}}}}}
\newcommand{\Tl}{{\Tcirc{\mbox{\raisebox{0.18ex}{$\scriptstyle [l]$}}}}}
\newcommand{\Ttwo}{{\Tcirc{\mbox{\raisebox{0.18ex}{$\scriptstyle [2]$}}}}}
\newcommand{\Ck}{{\mathrm{G}_k}}
\newcommand{\Ckstar}{{\mathrm{G}^*_k}}
\newcommand{\Ckinv}{{\mathrm{G}_k^{-1}}}
\newcommand{\Cke}{{\mathrm{G}_{k+1}}}
\newcommand{\Ckz}{{\mathrm{G}_{k+2}}}
\newcommand{\Ckl}{{\mathrm{G}_{k+l}}}
\newcommand{\Cl}{{\mathrm{G}_l}}

\newcommand{\Pk}{{\operatorname{P}_k}}
\newcommand{\Pke}{{\operatorname{P}_{k+1}}}
\newcommand{\Pkz}{{\operatorname{P}_{k+2}}}
\newcommand{\Pkl}{{\operatorname{P}_{k+l}}}

\newcommand{\hk}{{\mathrm{h}_k}}

\newcommand{\imc}{{\mathrm{imc}}}

\newcommand{\etafct}{{\eta}}

\newcommand{\lv}{\operatorname{lv}}

\newcommand{\istarkal}{{\mathrm{i}^\star_{k,\al}}}
\newcommand{\istarkbe}{{\mathrm{i}^\star_{k,\be}}}

\newcommand{\rstar}{{r^\star}}
\newcommand{\lstar}{{l^\star}}
\newcommand{\hstar}{{h^\star}}

\newcommand{\zeroseq}{\operatorname{0-Seq}}

\newcommand{\alcirc}{{\al^\circ}}

\begin{document}

\title{Generalizing Goodstein's theorem and Cichon's independence proof}

\author{
Gunnar Wilken\footnote{ORCID iD: 0000-0002-4019-9320, E-mail: wilken@oist.jp}\\
Okinawa Institute of Science and Technology, Japan\\
}

\maketitle

\begin{abstract} We generalize Goodstein's theorem \cite{Goodstein} and Cichon's independence proof \cite{Cichon1983} to $\pioneonecanod$ using
results from Wilken \cite{W26}. The method is generalizable to stronger notation systems that provide or enable unique terms for ordinals and enjoy Bachmann property, 
also abstractly via the uniform concept of maximality quotients.
\end{abstract}

\section{Introduction} 

Let $\Esc$ be the familiar notation system for ordinals below $\epsn$, the proof-theoretic ordinal of Peano arithmetic $\PA$,  
built from $0$, $(\xi,\eta)\mapsto\xi+\eta$, and $\xi\mapsto\om^\xi$ using Cantor normal form, where $\om$ denotes the least infinite ordinal number. 
For $\al\in\Esc$ let $\mc(\al)$ be the maximum counter of multiples that occur in $\al$, that is, the maximum number of times the same summand is consecutively added 
in the notation of $\al$. For $k\in[2,\om)$ define the quotient $\Equok\subseteq\Esc$ by \[\Equokbig:=\{\al\in\Esc\mid\mc(\al)<k\}.\]
Defining the system $(\Esc,\cdot[\cdot])$ of fundamental sequences (limit approximations) for $\Esc$ by
\[0[k]:=0,\;\;1[k]:=0,\;\;(\xi+\eta)[k]:=\xi+\eta[k],\;\;\om^{\be+1}[k]:=\om^\be\cdot(k+1),\;\mbox{ and }\;\om^\la[k]:=\om^{\la[k]},\]
and the slow-growing hierarchy $\Ck$ for $k>0$ by\footnote{Induction on $\al$ shows that our definitions are equivalent to Cichon's in \cite{Cichon1983} for $k>0$, where 
Cichon's $\al[k]$ is equal to our $\al[k-1]$.} 
\begin{equation}\label{slowgrowingeq} \Ck(0):=0,\: \Ck(\al+1):=\Ck(\al)+1, \mbox{ and }\Ck(\la):=\Ck(\la[k-1])\mbox{ for }\la\in\Esc\cap\Lim,\end{equation} 
where $\Lim$ denotes the class of limit ordinals, as shown in \cite{Cichon1983} the restricted mapping \[\Ck:\Equokbig\to\om\]
is a bijection via transforming back and forth the bases $k$ and $\om$, that is, writing a natural number $N$ in hereditary base-$k$ representation using $0$, addition,
and exponentiation to base $k$, replacing the base $k$ by $\om$ we obtain $\Ckinv(N)$, the unique preimage of $N$ in $\Equok$. 
The image of $\Equok$ under $\Ck$ is the Mostowski collapse of $\Equok$ and equal to $\om$, i.e., the restriction of $\Ck$ to $\Equok$ is an order isomorphism.

Since $(\Equok)_{k\in[2,\om)}$ is $\subseteq$-increasing, base transformation 
from $k$ to $l$ for $2\le k<l<\om$ is characterized by \[N[k\mapsto l]=\Cl(N[k\mapsto\om]),\]
where $N<\om$ is assumed to be given in hereditary base-$k$ representation, so that $N[k\mapsto\om]\in\Equok\subseteq\Equol$. Thus, ordinals $\al\in\Esc$ can be
seen as direct limit representations of natural numbers in hereditary base-$k$ representation, provided that $\mc(\al)<k$. 

In this article we generalize this mechanism to the notation system $\T$ for ordinals below Takeuti ordinal, that is, to the strength of the fragment $\pioneonecanod$ of second
order number theory, which is the strongest of the ``big five'' theories in Simpson's book on reverse mathematics, using fundamental
sequences and machinery introduced in \cite{W26}. The quotients introduced specifically for $\T$ restrict the iterations involved in limit approximation in a way directly corresponding
to the $\Equok$, by the introduction of a measure called {\it iterative maximal coefficient}, $\imc(\al)$, of an ordinal $\al\in\T$.
Such quotients of restricted iteration, in short $\imc$-quotients, reveal deep insight into the notation system $\T$.

Adopting the preliminaries of \cite{W26}, subsections 2.1, 2.3, and 2.5, which are summarized here in subsection \ref{prelimsec} for reader's convenience, 
let $\T$ be the notation system $\Tt$ for $\tau=1$ introduced there and let $\Tcirc$ be the
subset of terms of countable cofinality. The system $(\Tcirc,\cdot[\cdot])$ of fundamental sequences with Bachmann property given in Definition 3.5 of \cite{W26}
is called a Buchholz system. Bachmann property is a nesting property of fundamental sequences, namely that for terms $\al$ and $\be$ such that $\al[\ze]<\be<\al$ 
where $\al$ is not a successor-multiple of a regular cardinal $\Omie$, we have $\al[\ze]\le\be[0]$, see Theorem \ref{Bachmannprop}.


We will be able to use the approach by Cichon \cite{Cichon1983}, transferring the definition of collapsing functions $\Ck$ as defined in (\ref{slowgrowingeq})  
and predecessor functions $\Pk$ for $k>0$, defined equivalently as in \cite{Cichon1983} by
\begin{equation}\label{predeq} \Pk(0):=0,\: \Pk(\al+1):=\al, \mbox{ and }\Pk(\la):=\Pk(\la[k-1])\mbox{ for }\la\in\Lim,\end{equation}
to the Buchholz system $(\Tcirc,\cdot[\cdot])$, more specifically to its countable initial segment. 
Note first of all that Cichon's crucial lemma also holds in our context, with the same proof by induction on $\al$:
\begin{lem}[cf.\ Lemma 2 of \cite{Cichon1983}]\label{tricklem} For $k\in[2,\om)$ and $\al\in\T\cap\Om_1$ the operations $\Ck$ and $\Pk$ commute, thus \[\Ck\Pk\al=\Pk\Ck\al,\]
where the latter is equal to $\Ck\al-1$.
\end{lem} 

We are going to introduce quotients $\Tquok$ and $\Tcircquok$ of $\T$ and $\Tcirc$, respectively, to base $k$ where $k\in[2,\om)$. The sets $\Tk:=\Tcircquok\cap\Om_1$, 
each of which is cofinal in $\T\cap\Om_1$ and which are $\subseteq$-increasing in $k$, will be the canonical analogues of the quotients $\Equok$ mentioned above. 
It will be seen in Corollary \ref{maincor} that the predecessor function $\Pk$ is actually the true predecessor function on the restricted domain $\Tk$.
 
Consider as an instructive basic example the tower $k_k$ of height $k$ of exponentiation to base $k$.
Clearly, $k_k[k\mapsto\om]\in\Equok$, so that the direct limit corresponding to $k_k$ in $\Esc$ is $\om_k$, the tower of height $k$ of exponentiation to base $\om$.
This does not hold in $\Tcircquok$, where the corresponding ideal object is $\epsn=\lim_{n<\om}\om_n$. 
The sets $\Tk$ therefore provide more refined hierarchies of 
ideal objects uniquely denoting natural numbers, and functions $\Ck$ will act as their enumeration functions, see Theorem \ref{maintheo}.  

These preparations allow us to quite canonically generalize Goodstein's theorem \cite{Goodstein} and Cichon's independence proof in \cite{Cichon1983} to obtain independence
of the theory $\pioneonecanod$. 
The argumentation, ``{\it Cichon's trick}'', is as follows: Given a starting base $k\in[2,\om)$ and a natural number $N\in\N$, 
by Theorem \ref{maintheo} we obtain a {\it unique} $\al\in\Tk$ such that \[N=\Ck(\al)=:N_1.\] 
The uniqueness of $\al$ such that $\Ck(\al)=N$ in $\Tk$ allows us to understand $\al$ as {\it the base-$k$ representation of $N$}, 
and gives rise to the base transformation $N[k\mapsto l]:=\Cl(\al)$, so that the Goodstein operation of incrementing the base of representation and subsequent
subtraction of $1$ can be defined by \[N_1[k\mapsto k+1]-1=\Cke(\al)-1=\Pke\Cke\al=\Cke\Pke\al=:N_2,\] 
where we used Lemma \ref{tricklem}. Iterating the procedure yields, noting that $\Pke\al\in\Tke$, \[N_2[k+1\mapsto k+2]-1=\Ckz(\Pke\al)-1=\Pkz\Ckz(\Pke\al)=\Ckz\Pkz(\Pke\al)=:N_3,\]
so that we obtain the generalized Goodstein sequence \[N_{l+1}=\Ckl(\Pkl\ldots\Pke\al).\]
Termination of the generalized Goodstein sequence is therefore expressed by \[\exists l\:N_{l+1}=0,\]
and noting that $\Ck\al=0$ if and only if $\al=0$, we may reformulate the generalized Goodstein principle as follows:
\begin{equation}\label{reformulation}\forall k\in[2,\om)\:\forall\al\in\Tk\:\exists l\:\:\Pkl\ldots\Pke\al=0.\end{equation}
Thus, transfinite induction up to Takeuti ordinal proves this generalized Goodstein principle.
Independence then follows, once we show that the function $\hk:\T\cap\Om_1\to\N$ defined by 
\begin{equation}\label{hkdefieq}\hk(\al):=k + \min\{l\mid\Pkl\ldots\Pke\al=0\}\end{equation}
can easily be expressed in terms of the Hardy function $H_\al$, as in \cite{Cichon1983} for the original Goodstein principle. Defining Hardy hierarchy along Takeuti ordinal 
(i.e.\ the initial segment $\T\cap\Om_1$) in the same way as in \cite{Cichon1983} by 
\begin{equation}\label{Hardydefi} H_0(x):=x,\:H_{\al+1}(x):=H_\al(x+1), \mbox{ and }H_\la(x):=H_{\la[x]}(x),\end{equation}
we obtain by straightforward $<$-induction on $\al$:
\begin{equation}\label{hkHardyeq} \hk(\al)=H_\al(k).\end{equation}

The quotient method is generalizable to stronger notation systems that provide or enable unique terms for ordinals and enjoy Bachmann property.
One way is to develop quotients of restricted iteration for a stronger notation system, again leading to deeper insight into the specific notation system.

Another way yields an abstract, completely uniform generalization of Goodstein's theorem and Cichon's independence proof to
a system $(\tau,\cdot[\cdot])$ where $\tau$ is the proof-theoretic ordinal of a theory $\mathcal{T}$, which we may assume to be $\Pi^1_1$-sound and for which we may assume that a
$\Pi^0_2$-ordinal analysis is available, and $\cdot[\cdot]$ is a system of fundamental sequences that enjoys
Bachmann property, see Remark \ref{uniformityrmk}. 
In such context we use {\it maximality quotients}, see subsection \ref{quotientssubsec} for details. Independence of the resulting abstract Goodstein principle
then also follows from the unprovability of totality of the Hardy function $H_\tau$ in $\mathcal{T}$, which is assumed to hold.
This more general, uniform result has emerged from earlier versions of the present article and has been elaborated in more detail in \cite{W}. 

\section{Ordinal algebra}

\subsection{Preliminaries}\label{prelimsec}

Besides denoting the class of limit ordinals as $\Lim$, we denote the class of additive principal numbers as $\Hz$, $\Hz^{>1}:=\Hz\setminus\{1\}$, 
and the class of (non-zero) ordinals closed under $\om$-exponentiation, also called $\eps$-numbers, by $\Ez$.
 
For terms $\al$ we use the following abbreviations. Writing $\al=_\NF\xi+\eta$ means that $\eta\in\Hz$ and $\xi\ge\eta$ is minimal such that $\al=\xi+\eta$,
$\al=_\ANF\xi_1+\ldots+\xi_k$ means that $\xi_1,\ldots,\xi_k\in\Hz$ and $\xi_1\ge\ldots\ge\xi_k$, and to indicate Cantor normal form representation of an ordinal $\al$, 
we write $\al=_\CNF\om^{\ze_1}+\ldots+\om^{\ze_k}$ where $\ze_1\ge\ldots\ge\ze_k$.
For the result of $l$-fold addition of the same term $\eta$ to a term $\xi$, where $l<\om$, 
we will sometimes use the shorthand $\xi+\eta\cdot l$. We further define $\sumend(0):=0$ and $\sumend(\al):=\xi_k$ if $\al=_\ANF\xi_1+\ldots+\xi_k$ where $k\ge1$.

Setting $\Om_0:=\om$ and $\Omie:=\aleph_{i+1}$, terms in $\T$ are built up from $0$, $1$, ordinal addition, and functions $\thti$ for each $i<\om$, where $\thti(0)=\Omi$ 
and arguments of $\thti$ are restricted to terms below $\Omiz$. The functions $\thti$ are natural extensions of the fixed point-free Veblen functions and hence injective, 
with images contained in the intervals $[\Omi,\Omie)$, respectively for $i<\om$. This system $\T$ of $\thti$-functions therefore provides a system of stepwise collapsing functions, 
where defining $\theta_0:=\thtnod(0)=\om$ and $\theta_{n+1}:=\thtnod(\ldots(\tht_{n+1}(0)))$, the ordinal $\theta_{n+1}$ is well-known to be the proof-theoretic ordinal of the 
system $\IDn$ of $n$-times iterated inductive definitions where $\IDnod:=\PA$. 

Compared to the setup in \cite{W26}, defining $\Om_0:=\om$ instead of $\Om_0=1$ 
and consequently $\thtnod(i):=\om^{i+1}$ instead of $\thtnod(i)=\om^i$ for $i<\om$ is the only modification, which is needed for technical smoothness. 
In this context we define $\Om_0^\prime:=1$ and $\Omipr:=\Omi$ for $i>0$.

As in subsection 2.1 of \cite{W26}, slightly abusing notation we may consider notation systems 
$\Tomie$ to be systems relativized to the initial segment $\Omie$ of ordinals and built up over 
$\Omie=\thtie(0), \Om_{i+2}=\tht_{i+2}(0),\ldots$, i.e.\ without renaming the indices of $\tht$-functions. In this terminology $\T=\Tomnod$, and the sequence 
$(\Tomi)_{i<\om}$ can be seen as an increasing sequence with larger and larger initial segments of ordinals serving as parameters.

The operation $\cdot^\stari$ searches a $\T$-term for its $\thti$-subterm of largest ordinal value, but under the restriction to treat $\thtj$-subterms 
for $j<i$ as atomic. If such largest $\thti$-subterm does not exist, $\cdot^\stari$ returns $0$.  For $\al=\thti(\xi)$ and $\be=\thti(\eta)$ (where $\xi$ and $\eta$ are
any elements of the domain of $\thti$) we have \[\al<\be\Longleftrightarrow(\xi<\eta\mbox{ and }\xi^\stari<\be)\mbox{ or }\al\le\eta^\stari.\]

As in \cite{W26} the arguments of $\thti$-terms $\thti(\xi)$ are often written in a form $\thti(\De+\eta)$, i.e.\ $\xi=\De+\eta$. This by convention always indicates that
$\Omie\mid\De$ (possibly $\De=0$) and $\eta<\Omie$. 
Informally, we sometimes call the multiple $\De$ of $\Omie$ the fixed point level of $\thti(\xi)$, and the 
predicate $F_i(\De,\eta)$ means that $\eta=\sup^+_{\si<\eta}\thti(\De+\si)$, which according to Proposition 2.6 of \cite{W26} holds if and only if $\eta$ is of a 
form $\thti(\Ga+\rho)$ where $\Ga>\De$ and $\eta>\De^\stari$.

Any $\al\in\T\cap\Hz^{>1}$ is of a form $\thti(\De+\eta)$ for unique $i<\om$, $\De\in\T\cap\Omiz$ with $\Omie\mid\De$, and $\eta<\Omie$. 
Let $\lv(\al):=i$ denote the {\it level of $\al$}, and for $j\le i$ let $\Par_j(\al)$ denote the set of $\tht$-subterms of level strictly less than $j$ 
(so that $\Par_0(\al)=\emptyset$). 
We will sometimes write $\al\cdot\om$ as a shorthand for the additive principal successor of $\al=\thti(\De+\eta)$ in $\T$, 
instead of $\thti(\al)$ in case of $\De>0$ and $\thti(\eta+1)$ if $\De=0$.

For terms $\al$ of a form $\thti(\De+\eta)$ the $\Omi$-localization, cf.\ subsection 2.3 of \cite{W26}, is a sequence $\Omi=\al_0,\ldots,\al_m=\al$ ($m\ge0$), 
in which for $\al>\Omi$ the terms $\al_i=\thti(\De_i+\eta_i)$, $i=1,\ldots,m$, provide a strictly increasing sequence of $\thti$-subterms of $\al$ that are not in the scope of any 
$\thtj$-function with $j<i$ (i.e.\ that are not collapsed), such that the sequence of fixed point levels $\De_1,\ldots,\De_m$ is strictly {\it decreasing} and of maximal length,
and where the $\eta_i$ are maximal. 
Localization therefore approximates the ordinal $\al$ from below in terms of decreasing fixed point levels, and key properties are as follows.

\begin{prop}[cf.\ 2.21 of \cite{W26}]\label{localipic}
Let $\al=\thti(\De+\eta)\in\T$, $\al>\Omi$, with $\Omi$-localization $(\Omi=\alnod,\ldots,\alm=\al)$ where $\al_j=\thti(\De_j+\eta_j)$ for $j=1,\ldots,m$. Then 
\begin{enumerate} 
\item For $j<m$ and any $\be=\thti(\Ga+\rho)\in(\al_j,\al_{j+1})$ we have 
$\Ga+\rho<\De_{j+1}+\eta_{j+1}$.
\item $(\De_j)_{1\le j\le m}$ forms a strictly descending sequence of multiples of $\Omie$.
\item For $j<m$ the sequence $(\al_0,\ldots,\al_j)$ is the $\Omi$-localization of $\al_j$.
\end{enumerate} 
We have the following {\bf guiding picture} for localizations, setting $\al_j^+:=\thti(\De_j+\eta_j+1)$:
\[\Omi<\al_1<\ldots<\al_m=\al<\alplus=\al_m^+<\ldots<\al_1^+.\]
\end{prop}

Restricting ordinal addition occuring in $\T$ to summation in additive normal form, the terms in $\T$ {\it uniquely} identify ordinal numbers below Takeuti ordinal.
The term algebra developed in \cite{W26} and here therefore can be understood as ``algebra'' of ordinals. Further references are given in \cite{W26}.
\subsection{Fundamental sequences}
For reader's convenience we recall from \cite{W26} the following crucial definition of a characteristic function regarding uncountable moduli and then 
provide the concrete system of fundamental sequences adapted from Definition 3.5 of \cite{W26} to fit with our modified definition of $\thtnod(i)$ for $i<\om$.   
\begin{defi}[3.1 and 3.3 of \cite{W26}]\label{chidefi}
We define characteristic functions $\chiomie:\Tomie\to\{0,1\}$ for $i<\om$ by recursion on the build-up of $\Tomie$:
\begin{enumerate}
\item $\chiomie(\al):=\left\{\begin{array}{cl} 
                           0&\mbox{ if } \al<\Omie\\
                           1&\mbox{ if } \al=\Omie,
           \end{array}\right.$ 
\item $\chiomie(\al):=\chiomie(\eta)$ if $\al=_\NF\xi+\eta$,
\item $\chiomie(\al):=\left\{\begin{array}{ll}
          \chiomie(\De) & \mbox{ if } \eta\not\in\Lim\mbox{ or }F_j(\De,\eta)\\
          \chiomie(\eta)& \mbox{ otherwise,}
      \end{array}\right.$\\[1mm]
if $\al=\thtj(\De+\eta)>\Omie$ and hence $j\ge i+1$.
\end{enumerate}
For $\al\in\T$ we define $\domf(\al):=i+1$ if $\chiomie(\al)=1$ for some $i<\om$ and $\domf(\al):=0$ otherwise, so that \[\Tcirc=\{\al\in\T\mid\domf(\al)=0\}\]
is the subset of $\T$ containing the terms of ordinals of (at most) countable cofinality.
\end{defi}

Note that we obtain a partitioning of $\T$ in terms of cofinality through the preimages $\domf^{-1}(n)$, $n<\om$, cf.\ Lemma 3.4 of \cite{W26}.
We separately define the {\it support term} $\ual$ of a $\tht$-term $\al$ first, as follows.
\begin{defi}[cf.\  3.5 of \cite{W26}]\label{supporttermdefi} Let $\al\in\T$ be of the form $\al=\thti(\De+\eta)$ for some $i<\om$, and denote the 
$\Omi$-localization of $\al$ by $\Omi=\al_0,\ldots,\alm=\al$. The support term $\ual$ for $\al$ is defined by
\[\ual:=\left\{\begin{array}{cl}
            \almmin & \mbox{ if either } F_i(\De,\eta)\mbox{, or: } \eta=0\mbox{ and }\De[0]^{\star_i}<\almmin=\De^{\star_i}
            \mbox{ where }m>1\\
            \thti(\De+\etapr) & \mbox{ if } \eta=\etapr+1\\
            1 & \mbox{ if } \al=\om\\
            0 & \mbox{ otherwise.}
     \end{array}\right.\] 
\end{defi}
\begin{defi}[cf.\ 3.5 of \cite{W26}]\label{bsystemdefi}
Let $\al\in\T$. By recursion on the build-up of $\al$ we define the function \[\al[\cdot]:\aleph_d\to\T^{\Om_d}\] where $d:=\domf(\al)$.
Let $\ze$ range over $\aleph_d$.
\begin{enumerate}
\item $0[\ze]:=0$ and $1[\ze]:=0$.
\item $\al[\ze]:=\xi+\eta[\ze]$ if $\al=_\NF\xi+\eta$.
\item For $\al=\thti(\De+\eta)$ where $i<\om$, noting that $d\le i$, the definition then proceeds as follows.
\begin{enumerate}
\item[3.1.] If $\eta\in\Lim$ and $\neg F_i(\De,\eta)$, that is, $\eta\in\Lim\cap\sup_{\si<\eta}\thti(\De+\si)$, we have $d=d(\eta)$ 
and define 
   \[\al[\ze]:=\left\{\begin{array}{cl}
   \thtnod(\ze)\ & \mbox{ if } \al=\thtnod(\thtnod(0))=\om^\om\\
   \thti(\De+\eta[\ze]) & \mbox{ otherwise.}
   \end{array}\right.\]
\item[3.2.] If otherwise $\eta\not\in\Lim$ or $F_i(\De,\eta)$, we distinguish between the following 3 subcases.
\begin{enumerate}
\item[3.2.1.] If $\De=0$, define \[\al[\ze]:=\left\{\begin{array}{cl} \ual\cdot(\ze+1) & \mbox{ if } d=0,\\
            \ze & \mbox{ if } d>0.
       \end{array}\right.\]
\item[3.2.2.] $\chiomie(\De)=1$. This implies that $d=0$, and we define recursively in $n<\om$
\[\al[0]:=\thti(\De[\ual])\quad\mbox{ and }\quad\al[n+1]:=\thti(\De[\al[n]]).\]
\item[3.2.3.] Otherwise. Then $d=d(\De)$ and \[\al[\ze]:=\thti(\De[\ze]+\ual).\]
\end{enumerate}
\end{enumerate} 
\end{enumerate}
\end{defi}

The following two lemmas from \cite{W26} are included for the reader's convenience and address the main (technical) properties of fundamental sequences for ordinals in $\T$ and  
$\Omi$-localizations of elements of fundamental sequences. 

\begin{lem}[cf.\ Lemma 3.10 of \cite{W26}]\label{bracketsmainlem} 
In the setting of Definition \ref{bsystemdefi}, let $\al\in\Tt\setminus\tau$ be a limit ordinal, and let $\ze$ range over $\aleph_d$ 
where $d:={\domf(\al)}$.
\begin{enumerate}
\item\label{ualcontrolclaim} For $\al=\thtj(\De+\eta)$ such that $\ual>0$ and $\chiomje(\De)=0$ we have \[\De[\ze]^\starj<\ual.\]
\item\label{fundseqclaim} The mapping $\ze\mapsto\al[\ze]$ is strictly increasing with proper supremum $\al$. 
\item\label{starcontrolclaim} The mapping $\ze\mapsto\al[\ze]^{\stark}$ is weakly increasing with upper bound $\al^{\stark}$
 for $d\le k$. 
\item\label{bracketparestimclaim} If $d=i+1>0$, we have 
\[\ze^\stark\le\al[\ze]^{\stark}\le\max\{\al^{\stark},\ze^{\stark}\}\quad\mbox{ and }\quad\al^{\stark}\le\al[\ze]^{\stark}\] 
for any $k\le i$.
\item\label{bracketlocalizationclaim} For $\al=\thtj(\De+\eta)>\Omj$, denote the $\Omj$-localization of $\al$ by 
$\Omj=\alnod,\ldots,\alm=\al$. For such $\al$ we have \[\almmin\le\al[0].\]
\item\label{stardomclaim} If $\al=\thtj(\De+\eta)$ and $i$ is such that $\domf\le i\le j$ and $\al[0]^\stari<\al^\stari$
where $\al^\stari>\Omi$, then we have $\domf(\al^\stari)=\domf$, and one of the following two cases applies.
\begin{enumerate}
\item $\al^\stari$ is of the form $\nu\cdot\om$ for some additive principal $\nu\ge\Omi$, where we have $d=0$,
$\al[n]^\stari=\nu$, and $\al^\stari[n]=\nu\cdot(n+1)$ for $n<\om$. 
\item Otherwise. Then $\al^\stari[\ze]\le\al[\ze]^\stari$ for all $\ze<\aleph_\domf$, 
and there exists $\ze_0<\aleph_\domf$ such that $\al^\stari[\ze]=\al[\ze]^\stari$ for all $\ze\in(\ze_0,\aleph_\domf)$.
\end{enumerate}
\item\label{sandwichclaim} For any $\be$ such that $\al[\ze]\le\be<\al$ we have \[\be^\stark\ge\al[\ze]^\stark\] 
for all $k$ if $\ze=0$, and 
for all $\ze<\aleph_\domf$ and all $k$ such that $k+1\ge\domf$.
\end{enumerate}
\end{lem}

We include a remark from \cite{W26} as a proper corollary since it is a deep property of the system $\cdot[\cdot]$ on $\T$ that will be used
in Lemmas \ref{zeroiterationlem}, \ref{supportbyzeroiterationlem}, and \ref{extendedTquoksubstlem}.
 
\begin{cor}[cf.\ Remark 3.8, part 6, of \cite{W26}]\label{deepsupporttermcor} 
Let $\al=\thti(\De+\eta)$ with $\Omi$-localization $(\Omi,\al_1,\ldots,\alm=\al)$.
Suppose that $\eta=0=\ual$ and $\De$ is a limit-multiple of $\Omie$ such that $\chiomie(\De)=0$. Then we observe:
\begin{enumerate}
\item The situation $\De[0]^\stari<\almmin<\De^\stari$ does not occur. 
\item If  $m>1$ and $\De[0]^\stari<\De^\stari$, it follows that $\almmin\le\De[0]^\stari$.
\end{enumerate}
\end{cor}
{\bf Proof.}
Part 1 follows from parts \ref{bracketlocalizationclaim} and \ref{stardomclaim} of Lemma \ref{bracketsmainlem}. Assuming that $\De[0]^\stari<\almmin<\De^\stari$,
it follows that $\De^\stari>\Omi$, and, as $\domf(\De)\le i<i+1$,  part \ref{stardomclaim} (b) of Lemma \ref{bracketsmainlem} applies and yields
$\almmin\le\De^\stari[0]\le\De[0]^\stari$, where we applied part  \ref{bracketlocalizationclaim} to $\De^\stari$. Contradiction.

For part 2, assuming additionally that $m>1$ and $\De[0]^\stari<\De^\stari$, we observe that $\almmin<\De^\stari$ as ``$=$'' is excluded by the assumption $\eta=0=\ual$.
But according to part 1, we can not have $\De[0]^\stari<\almmin<\De^\stari$.
\qed

\begin{lem}[cf.\ 3.12 of \cite{W26}]\label{localizationlem}
Let $\al=\thti(\De+\eta)>\Omi$ with $\Omi$-localization $(\Omi=\alnod,\ldots,\alm=\al)$, and let $\ze<\aleph_{\domf(\al)}$.
We display the $\Omi$-localization of $\al[\ze]$, distinguishing between the cases of definition. 
\begin{enumerate}
\item $\eta\in\Lim$ such that $F_i(\De,\eta)$ does not hold. Then the $\Omi$-localization of $\al[\ze]$ is 
$(\alnod,\ldots,\almmin,\al[\ze])$, unless $\ze=0$ and either $\al=\thtnod(\thtnod(0))$ or $\al=\thti(\Omk)$ with $0<k\le i$, where the $\Omi$-localization of
$\al[0]$ is $(\Omi)$.
\item $\eta\not\in\Lim$ or $F_i(\De,\eta)$ holds.
\begin{enumerate}
\item $\De=0$. Here $\al[\ze]\in\Hz$ if and only if $\ze=0$. Since $\al>\Omi$, the $\Omi$-localization of $\al[0]$ then is $(\alnod,\ldots,\almmin)$,
unless $\eta=\etapr+1$ for some $\etapr>0$, in which case it is $(\alnod,\ldots,\almmin,\ual)$ where $\ual=\thti(\etapr)$.
\item $\chiomie(\De)=1$. The $\Omi$-localization of $\al[\ze]$ then is $(\alnod,\ldots,\almmin,\ual,\al[\ze])$ if $\eta=\etapr+1$
for some $\etapr$ where $\ual=\thti(\De+\etapr)$, otherwise it is $(\alnod,\ldots,\almmin,\al[\ze])$, unless $\ze=0$, $\De=\Omie$, 
and $\eta=0$, in which case it simply is $(\Omi)$.
\item $\De>0$ and $\chiomie(\De)=0$. Then the $\Omi$-localization of $\al[\ze]$ is $(\alnod,\ldots,\almmin,\al[\ze])$, unless
$\eta=\etapr+1$ for some $\etapr$, in which case it is $(\alnod,\ldots,\almmin,\ual,\al[\ze])$ where $\ual=\thti(\De+\etapr)$.
\end{enumerate}
\end{enumerate}
\end{lem}

We now consider the iterated descent via fundamental sequences and derive useful lemmas.

\begin{defi}\label{iterateddescentdefi} 
\begin{enumerate}
\item Let $l<\om$ and $\al\in\T\cap\Om_1$. We define the iterated application of $\cdot[l]$ by 
\[\al[l]^0:=\al\;\mbox{ and }\;\al[l]^{j+1}:=(\al[l]^j)[l] \mbox{ for }l<\om.\]
For $k>0$ let $\istarkal<\om$ to be minimal such that $\al[k-1]^\istarkal=0$. 
\item For arbitrary $\al\in\T$ and suitable $\ze$ the iterated application of $\cdot[\ze]$ is defined likewise as far as this makes sense. 
\item We define in particular the sequence of iterated application of $\cdot[0]$ to $\al\in\T$ as \[\zeroseq(\al):=(\al[0]^h)_{h\le\hstar}\]
where $\hstar$ is minimal such that $\al[0]^\hstar=0$. 
\[G(\al):=\hstar\] is called the {\it G-norm} of $\al$, cf.\ Definition 5.4 of \cite{W26}.
\end{enumerate}
\end{defi}

\begin{lem}\label{zeroiterationlem} 
Let $\Ga=\thtje(\Xi+\rho)$ be such that $\chiomie(\Ga)=1$ and $\be\in\T\cap\Lim\cap[\Omi,\Omie)$. Then we have
\begin{enumerate}
\item $\domf(\Ga[\be])=\domf(\be)\le i$, and, in case of $i=0$, $\domf(\Ga[\ze+1])\le j+1$ for $\ze<\Om_1$.
\item If $\Ga$ is a limit-multiple of $\Omje$, then so are $\Ga[\be]$ and (in case of $i=0$) $\Ga[\ze+1]$ for $\ze<\Om_1$.
\item $\Ga[\be][\ze]=\Ga[\be[\ze]]$ for $\ze<\aleph_{\domf(\be)}$.
\item There is an $\lstar$ such that $\Ga[\be][0]^\lstar=\left\{\begin{array}{cl}
\Ga[\be_0+1] & \mbox{ if }i=0, \mbox{ where }\be=_\NF\be_0+\sumend(\be)\\
\Ga[0] & \mbox{ if }i>0.
\end{array}\right.$ 
\end{enumerate}  
\end{lem}
{\bf Proof.} We prove all claims apart from the last one simultaneously by induction on the build-up of $\Ga$. Note that we have $i\le j$. The case $\Ga=\Omie$ is trivial.
Note that the iterated application of $\cdot[0]$ to $\sumend(\be)$ produces ordinals in $\Hz\cap[\Omi,\Omie)$ until reaching $\Omi$, so that another application of $\cdot[0]$ finally
reaches $1$ in case of $i=0$ and $0$ otherwise. Thus the last claim eventually follows.
\\[2mm]
\noindent{\bf Case 1:} $\rho\in\Lim$ and $F_{j+1}(\Xi,\rho)$ does not hold. Then we have $\domf(\Ga)=\domf(\rho)=i+1$ and $\Ga[\be]=\thtje(\Xi+\rho[\be])$, 
where $\rho[\be]\in\Lim$ and $F_{j+1}(\Xi,\rho[\be])$ does not hold, as follows from Lemma \ref{localizationlem}. Thus $\domf(\Ga[\be])=\domf(\rho[\be])$, 
$\Ga[\be][\ze]=\thtje(\Xi+\rho[\be][\ze])$, and applying the i.h.\ to $\sumend(\rho)$ we obtain $\domf(\Ga[\be])=\domf(\be)\le i$ and $\Ga[\be][\ze]=\Ga[\be[\ze]]$. 
Clearly,  in this case $\Ga$ and $\Ga[\be]$, as well as $\Ga[\ze+1]=\thtje(\Xi+\rho[\ze+1])$ in case of $i=0$, are limit-multiples of $\Omje$. Note that if $i=0$, $\domf(\Ga[\ze+1])\le j+1$ holds trivially.
\\[2mm]
\noindent{\bf Case 2:} $\rho\not\in\Lim$ or $F_{j+1}(\Xi,\rho)$ holds. Then we have $\domf(\Ga)=\domf(\Xi)$, which implies that we must have $\Xi>0$ such that $\chiomjz(\Xi)=0$ as
$i+1<j+2$. Hence, $\Ga[\be]=\thtje(\Xi[\be]+\uGa)$. The i.h.\ applies to $\sumend(\Xi)$, and as $\Xi$ is a limit-multiple of $\Omjz$, so is $\Xi[\be]$. 
We further have $\domf(\Ga[\be])=\domf(\Xi[\be])=\domf(\be)\le i$ and $\Xi[\be][\ze]=\Xi[\be[\ze]]$. 
We claim that $\underline{\Ga[\be]}=\uGa$. This is clear if $\uGa>0$. Otherwise we have $\rho=\uGa=0$. Let $\Ga_{m-1}$ be the predecessor of $\Ga$ in its $\Omje$-localization.
According to Lemma \ref{localizationlem} $\Ga_{m-1}$ then also is the immediate predecessor of $\Ga[\be]$ in its $\Omje$-localization. Assume that $m>1$, as otherwise 
$\underline{\Ga[\be]}=0$ would follow immediately. Then the assumption \[\Xi[\be][0]^\starje=\Xi[\be[0]]^\starje<\Xi[\be]^\starje=\Ga_{m-1}\] would entail $\Xi[0]^\starje<\Ga_{m-1}<\Xi^\starje$, 
as $\uGa=0$ and $\cdot^\starje$ is weakly increasing in this context. 
As this is impossible according to Corollary \ref{deepsupporttermcor}, it follows that $\underline{\Ga[\be]}=\uGa$ and hence $\Ga[\be][\ze]=\Ga[\be[\ze]]$. 
The claims regarding $\Ga[\ze+1]$ in case of $i=0$ follow as well, as in the case $\domf(\Xi[\ze+1])=j+2$ we have $\domf(\Ga[\ze+1])=0$, concluding the proof.
\qed

The following lemma precisely describes the pattern of descent via iteration of $\cdot[0]$. This will turn out to be crucial for our understanding of quotients.

\begin{lem}\label{supportbyzeroiterationlem}
Let $\al=\thti(\De+\eta)\in\T$ with $\Omi$-localization $(\Omi,\al_1,\ldots,\al_m=\al)$. 
\begin{enumerate}
\item There exists $h$ such that \[\ual=\al[0]^h.\]
\item The sequence $\zeroseq(\al)=(\al[0]^h)_{h\le\hstar}$ passes in reverse order through the $\Omi$-localization of $\al$, which therefore is a reversed subsequence.
$\zeroseq(\al)$ consists of ordinals in $\Hz^{\ge\Omi}\cup\{0,1\}$, which in case of $i=0$ always ends with elements $\om,1,0$, and in case of $i>0$ ends with elements $\Omi,0$.
\item Let $(\eta[0]^l)_{l\le\lpr}$ be equal to $\zeroseq(\eta)$, unless $m>1$ and the ordinal $\almmin$ is attained (this happens if and only if $\eta\ge\almmin>\Omi$), 
in which case we truncate the sequence  $\zeroseq(\eta)$ after $\almmin$. 
Then \[(\thti(\De+\eta[0]^l))_{l\le\lpr}\] is a subsequence of $\zeroseq(\al)$ and passes through ordinals $\thti(\De+\etapr)$ for initial sums $\etapr\ge\eta[0]^\lpr$ of $\eta$. 
Hence for $h=\lpr$ we have 
\[\al[0]^h=\left\{\begin{array}{cl} 
\thti(\De+\almmin) & \mbox{ if } \eta\ge\almmin>\Omi\\[2mm] 
\thti(\De) & \mbox{  otherwise.}\end{array}\right.\]
In the base case $\De=0$ we therefore obtain $\al[0]^{h+1}=\almmin$ if $\eta\ge\almmin>\Omi$, while otherwise $\al[0]^{h+1}=0$ if $i>0$ and $\al[0]^{h+1}=1$ if $i=0$.
\item If $\chiomie(\De)=1$, clearly $\thti(\De)[0]=\thti(\De[0])$ in the situation $\eta=0$, as then $\ual=0$. Suppose now that $F_i(\De,\eta)$ holds, hence $\eta=\almmin$ where $m>1$. 
We then have 
\[\al[0] = \thti(\De[\almmin]).\]
\begin{enumerate}
\item If $\De=_\NF\De[0]+\Omie$, i.e.\ $\De$ is a successor-multiple of $\Omie$, then $\al[0]=\thti(\De[0]+\almmin)$.
\item If $\De$ is a limit-multiple of $\Omie$, we have 
\[\al[0]^{2+h} = \thti(\De\left[\almmin[0]^{1+h}\right]+\almmin)\]
for $h\le\hstar$ where $\hstar$ is such that $\almmin[0]^{1+\hstar}=\left\{\begin{array}{cl} 1 & \mbox{ if } i=0\\ 0 & \mbox{  if }i>0.\end{array}\right.$  
\begin{itemize}
\item If $i>0$, we reach $\al[0]^{2+\hstar}=\thtnod(\De[0]+\almmin)$.
\item If $i=0$, we consider \[\alpr:=\al[0]^{2+\hstar}=\thtnod(\De[1]+\almmin).\] 
There exists an $h$ such that $\alpr[0]^h=\thti(\De[0]+\almmin)$.
\end{itemize}
\end{enumerate}
\item If $\De>0$ such that $\chiomie(\De)=0$, we consider cases regarding $\eta$.
\begin{enumerate}
\item If $F_i(\De,\eta)$ holds, hence $\eta=\almmin$ where $m>1$, we have $\al[0]=\thti(\De[0]+\almmin)$.
\item If $\eta=0$, let $\zeroseq(\De)=(\De[0]^h)_{h\le\hstar}$ and let $\hpr$ be minimal such that $(\De[0]^{\hpr+1})^\stari<\almmin$ where $m>1$, if that exists, and $\hpr:=\hstar$ otherwise.  
We then have
\[\al[0]^h=\thti(\De[0]^h)\]
for $h\le\hpr$, and in case of $\hpr<\hstar$
\[\al[0]^{\hpr+1}=\thti(\De[0]^{\hpr+1}+\almmin),\]
while otherwise $\al[0]^\hstar=\thti(\Omie)[0]=\Omi$.
\end{enumerate} 
\end{enumerate} 
\end{lem}
{\bf Proof.} 
\\[2mm]
{\bf Ad 1.} We argue by induction on $\al$.
Since the iterated application of $\cdot[0]$ to $\al$ yields a strictly descending sequence, the claim is trivial 
if $\ual=0$. Suppose $\ual>0$ from now. We either have $\eta\not\in\Lim$ or $F_i(\De,\eta)$ holds, hence
\[\al[0]=\left\{\begin{array}{cl}
\ual & \mbox{ if }\De=0\\[2mm]
\thti(\De[\ual]) & \mbox{ if }\chiomie{\De}=1\\[2mm]
\thti(\De[0]+\ual) & \mbox{ otherwise.}
\end{array}\right.\]
If $\De=0$ we are done, otherwise we consider the following cases.
\\[2mm]
{\bf Case 1:} $\ual=\almmin$. Then $\almmin\le\al[0]$ by part 5 of Lemma \ref{bracketsmainlem}, and according to
Lemma \ref{localizationlem} $\almmin$ is the immediate predecessor of $\al[0]$ in its $\Omi$-localization. Iteration of $\cdot[0]$ eventually reaches $\almmin=\ual$.
\\[2mm]
{\bf Case 2:} $\ual=\thti(\De+\etapr)$  where $\eta=\etapr+1$. Then according to Lemma \ref{localizationlem} $\ual$ is the immediate predecessor of $\al[0]$ in its $\Omi$-localization.
We claim that \[\underline{\al[0]}=\ual.\] 
{\bf 2.1:} $\chiomie(\De)=0$. Then we have $F_i(\De[0],\ual)$ since $\De[0]^\stari<\ual$, cf.\ part 1 of Lemma \ref{bracketsmainlem}. 
\\[2mm]
{\bf 2.2:} $\chiomie(\De)=1$. Then part 4 of Lemma \ref{bracketsmainlem} yields $\De[\ual]^\stari=\ual$ since $\De^\stari<\ual=\thti(\De+\etapr)$. 
\\[2mm]
{\bf 2.2.1:} $\De=_\NF\De[0]+\Omie$, i.e.\ $\De$ is a successor-multiple of $\Omie$. Then we have $\al[0]=\thti(\De[0]+\ual)$. 
\\[2mm]
{\bf 2.2.2:} Otherwise, that is, $\De$ is a limit-multiple of $\Omie$. Then by part 4 of Lemma \ref{bracketsmainlem} and Lemma \ref{zeroiterationlem} the ordinal $\De[\ual]$ is a 
limit-multiple of $\Omie$ with $\chiomie(\De[\ual])=0$ and
\[\De[\ual][0]^\stari=\De[\ual[0]]^\stari\le\max\{\De^\stari,\ual[0]\}<\ual=\De[\ual]^\stari,\] hence $\underline{\al[0]}=\ual$. 
\\[2mm]
We may now apply the i.h.\ to $\al[0]$ to obtain the claim.
\\[2mm]
{\bf Ad 2.} This follows from the definition of $\cdot[0]$, part 5 of Lemma \ref{bracketsmainlem}, and Lemma \ref{localizationlem}. 
\\[2mm]
{\bf Ad 3.} This is a consequence of the involved definitions, using part 1.
\\[2mm]
{\bf Ad 4.} This follows from Lemma \ref{zeroiterationlem}. For part b in the situation $i=0$, note that by Lemma \ref{zeroiterationlem} $\De[1]$ is a limit-multiple of $\Om_1$ 
such that $\domf(\De[1])\le1$. We argue by induction on the fixed-point level of $\alpr[0]^h$ while iterating the application of $\cdot[0]$. Due to Bachmann property we have
$\De[0]\le\De[1][0]^h$ as long as $\De[0]<\De[1][0]^\hpr$ for $\hpr<h$. This clearly also holds for $\hpr$ such that $\De[1][0]^\hpr=\Depr+\Om_1$ for some $\Depr$. Note that
``detours'' may occur when $\chiome(\De[1][0]^h)=1$ with $\De[1][0]^h$ a limit-multiple of $\Om_1$. For the least such $h$ we have $\alpr[0]^{h+1}=\thtnod(\De[1][0]^h[\almmin])$,
and by Lemma \ref{zeroiterationlem} further iteration of $\cdot[0]$ eventually yields $\thtnod(\De[1][0]^h[1]+\almmin)$. This pattern may repeat until this process eventually reaches an
ordinal $\alpr[0]^h$ with fixed-point level $\De[0]$, while the support term remains $\almmin$ throughout. 
\\[2mm]
{\bf Ad 5.} This follows from the definitions using Corollary \ref{deepsupporttermcor}. 
\qed

The above lemma has the following important corollary which clarifies the relationship between a $\tht$-term and its argument regarding the $G$-norm.

\begin{cor}\label{Gnormcor}
Let $\al=\thti(\De+\eta)\in\T$ with $\Omi$-localization $(\Omi,\al_1,\ldots,\al_m=\al)$. 
For any initial sum $\xi$ of $\De+\eta$ the sequence $\zeroseq(\al)$ passes either through $\thti(\xi)$, namely if $m\le1$ or $\xi^\stari\ge\almmin$, or through $\thti(\xi+\almmin)$ in
case of $m>1$ and $\xi^\stari<\almmin$. Moreover, we have \[G(\De+\eta)<G(\al).\]
\end{cor}

The above observations give rise to the following technical definition of a ``socle'' term. Consider a term $\al=\thti(\xi)\in\Tquok$ where $\xi=_\ANF\xi_1+\ldots+\xi_n\in\T\cap\Omiz$.
$\Tquok$ is not closed under ``initial sum $\tht$-terms'', i.e., we do not in general have $\thti(\xi_1+\ldots+\xi_l)\in\Tquok$ for $l<n$.   
Note that $\alcirc$ as defined below satisfies $\alcirc<\al$ and is an element of $\zeroseq(\al)$ by Corollary \ref{Gnormcor}, showing that $\Tquok$ is closed under the operation $\al\mapsto\al^\circ$.

\begin{defi}\label{socletermdefi}
Let $\al=\thti(\De+\eta)\in\T$ with $\Omi$-localization $(\Omi,\al_1,\ldots,\al_m=\al)$. We define the {\it socle} $\alcirc$ of $\al$ by the following case distinction:
\begin{enumerate}
\item If $\De=0$ and either $\eta=0$ or $F_i(0,\eta)$ holds, we set $\alcirc:=\eta$.
\item If $\eta=_\NF\etapr+\sumend(\eta)$ such that either $\eta=1$ or $\etapr>0$, i.e.\ $\eta$ is a successor or additively decomposable, we set $\alcirc:=\thti(\De+\etapr)$.
\item If $\eta\in\Hz^{>1}$ such that $F_i(\De,\eta)$ does not hold, then defining
$\uual:=\left\{\begin{array}{cl} 
\almmin & \mbox{ if }\De^\stari<\almmin<\eta\mbox{ with }m>1\\
0 & \mbox{ otherwise,}
\end{array}\right.$
we set $\alcirc:=\thti(\De+\uual)$.
\item Otherwise, we have $\De>0$, and either $\eta=0$ or $F_i(\De,\eta)$ holds. Then, writing $\De=_\NF\Depr+\sumend(\De)$, we set $\alcirc:=\thti(\Depr+\uual)$ where
$\uual:=\left\{\begin{array}{cl} 
0 & \mbox{ if }m=1\mbox{ or }\almmin\le(\Depr)^\stari\\
\almmin & \mbox{ otherwise.}
\end{array}\right.$
\end{enumerate}
\end{defi}

\subsection{Inversion}
By means of the following inversion lemma we will be able to recover $\be$ from terms $\al\in\T$ of a form $\be[n]$ where $n\in[1,\om)$ and $\be\in\Tcirc\cap\Lim$, 
as used in Definition \ref{quotientdefi}. As in Definition 3.5 of \cite{W26} we need to formulate the inversion lemma more generally, so as to cover terms of uncountable cofinality.
The inversion is formulated in a way that will cover the cases needed to define the quotients of Definition \ref{quotientdefi}, rather than covering all possible cases.
We first introduce a partial auxiliary function that will simplify the formulation of the subsequent lemma, cf.\ the support term $\ual$ in Definition \ref{supporttermdefi}.

\begin{defi} For $i<\om$, $\De\in\T$ such that $\Omie\mid\De<\Omiz$, and $\rho\in\T$ we define the partial function $\etafct(i,\De,\rho)$ as follows.
\[\etafct(i,\De,\rho):=\left\{\begin{array}{cl}
           \rho & \mbox{ if } F_i(\De,\rho) \mbox{ holds,} \\[2mm]
           0 & \mbox{ if either } \rho=0,\; \rho=1 \mbox{ and }i=\De=0, \mbox{ or: } \De[0]^\stari<\al_{m-1}=\De^\stari=\rho\\[1mm] 
              & \mbox{ where } \al_1,\ldots,\al_m \mbox{ is the $\Omi$-localization of } \thti(\De) \mbox{ and } m>1, \\[2mm]
           \nu+1 & \mbox{ if } \rho=\thti(\De+\nu) \mbox{ for some }\nu\in\T\cap\Omie, \\[2mm]
           \mbox{ undefined } & \mbox{ otherwise.}
      \end{array}\right.\]
\end{defi}

\begin{rmk}
In the cases where $\eta:=\etafct(i,\De,\rho)$ is defined, for $\al:=\thti(\De+\eta)$ we then have $\ual=\rho$, and either $\eta\not\in\Lim$ or $F_i(\De,\eta)$ holds.
\end{rmk}

\begin{lem}[Inversion Lemma]\label{inversionlem} Let $\al\in\T$. $\al$ is of a form $\al=\be[\ze]$ where 
\begin{enumerate}
\item $\ze\in[1,\om)$ and $\be\in\Tcirc\cap\Lim$, or 
\item $\ze\in[\Omi,\Omie)\cap\Hz$ for some $i<\om$ and $\be\in\T$ such that $\chiomie(\be)=1,$
\end{enumerate}
if and only if one of the following cases applies:
\begin{enumerate}
\item $\al=\ze\in[\Omi,\Omie)\cap\Hz$ and $\be=\Omie=\thtie(0)$.
\item $\al=\eta\cdot(n+1)$ for some $\eta\in\T\cap\Hz$ and $\ze=n\in[1,\om)$. Then $\be:=\eta\cdot\om$, and we have $\ube=\eta$ and $\al=\be[\ze]$.
\item $\al=\xi_1+\ldots+\xi_k+\rho$ where $\xi_1,\ldots,\xi_k\in\Hz$ is weakly decreasing, $k\ge1$, and $\rho=\eta[\ze]$ according to the lemma's conditions with
$\eta\in(1,\xi_k]\cap\Hz$ so that $\rho<\xi_k$. 
Then setting $\be:=\xi_1+\ldots+\xi_k+\eta$ we have $\al=\be[\ze]$.
\item $\al=\thtj(\Ga+\rho)$ and one of the following subcases applies:
\begin{enumerate}
\item Either $j=0$, $\Ga=0$, and $\rho=\ze\in[1,\om)$, so that $\al=\om^{1+\ze}=\be[\ze]$ for $\be:=\thtnod(\thtnod(0))=\om^\om$, or otherwise
$\rho$ is of a form $\rho=\eta[\ze]$ according to the lemma's conditions where $\eta<\Omje$ and $F_j(\Ga,\eta)$ does not hold. Then setting $\De:=\Ga$ and
$\be:=\thtj(\De+\eta)$ we have $\al=\be[\ze]=\thtj(\De+\eta[\ze])$.
\item Setting $\xi_1:=\Ga+\rho$, check whether there is a (shortest) sequence $\xi_1,\ldots,\xi_{m+1}$ (where $m\ge1$) 
that determines a term $\De$ with $\Omje\mid\De<\Omjz$ and $\chiomje(\De)=1$ in the first step, such that
\begin{enumerate}
\item $\xi_k$ is of a form $\De[\thtj(\xi_{k+1})]$ which is according to the lemma's conditions (here according to condition 2 with i=j) 
for $k=1,\dots,m$, and
\item $\xi_{m+1}$ is of a form $\De[\nu]$ where $\eta:=\eta(j,\De,\nu)$ is defined, so that $\be:=\thtj(\De+\eta)$, $\ube=\nu$, 
and $\al=\be[m]$. 
\end{enumerate}
This case then applies if  $\ze=m\ge1$.
\item $\Ga$ is of a form $\De[\ze]$ according to the lemma's conditions where $\Omje\mid\De<\Omjz$ and $\chiomje(\De)=0$, so that $\al=\thtj(\De[\ze]+\rho)$, 
and $\eta:=\eta(j,\De,\rho)$ is defined, so that $\be=\thtj(\De+\eta)$, $\ube=\rho$, and $\al=\be[\ze]$. 
\end{enumerate} 
\end{enumerate}
\end{lem}
{\bf Proof.} Correctness follows by induction on $\al\in\T$: If one of cases $1$ - $4$ holds with $\be$ and $\ze$, then $\al=\be[\ze]$ matches either condition $1$ or $2$.
The reverse direction, completeness, follows by induction on $\be\in\T$, showing that for any $\be\in\T$ and $\ze$ according to either condition $1$ or $2$, the ordinal
$\al:=\be[\ze]$ satisfies one of cases $1$ - $4$.  
\qed

\begin{rmk} 
Considering the example $\epsn$ as in the introduction, $\epsn[1]=\thtnod(\thtnod(0))=\om^\om$ is the first iterative repetition, while $\epsn[0]=\om=\thtnod(0)$, and 
$\om^\om[1]=\om^2=\thtnod(1)$. 
\end{rmk}

\subsection{Bachmann property}\label{Bachmannsubsec}
The following theorem states that terms in $\T$ satisfy Bachmann property, which is crucial for our approach toward generalizing Goodstein's theorem.
 
\begin{theo}[Bachmann property, Lemma 4.1 and Theorem 4.2 of \cite{W26}]\label{Bachmannprop}
For $\al, \be\in\T$, where $\al$ is not a successor-multiple of a regular cardinal, and for which $\al[\ze]<\be<\al$ holds for some $\ze<\aleph_{\domf(\al)}$, we have \[\al[\ze]\le\be[0].\] 
\end{theo}

The following corollary states a basic observation about fundamental sequences with Bachmann property. Note first that whenever $\al[\xi]=\al[\ze]$ where
$\al\in\T\cap\Lim$ and $\xi,\ze<\aleph_{\domf(\al)}$, then $\xi=\ze$ since fundamental sequences are strictly increasing. Note further that for any $\al\in\T\cap\Hz$
we have $\al\cdot\om[0]=\al$.

\begin{cor}\label{Bachmanncor} Let $\al,\be\in\T$ such that $\be<\al$ and $\al$ is not a successor-multiple of any $\Omje$. Let $\ze<\aleph_{\domf(\be)}$, $\xi<\aleph_{\domf(\al)}$,
and suppose that \[\al[\xi]=\be[\ze].\] Then it follows that $\ze=0$.
\end{cor}
{\bf Proof.} This is an immediate consequence of Bachmann property, Theorem \ref{Bachmannprop}, as assumptions imply that $\be[\ze]\le\be[0]$.
\qed

\begin{rmk}\label{nexusrmk} For ordinals $\al,\be\in\Tcirc\cap\Hz$ it naturally occurs that $\al[0]=\be[m]$ for some $m>0$. Due to Bachmann property this implies that $\al<\be$, 
and $\be$ is of greater fixed-point level than $\al$
unless $\be$ is of a form $\thti(\Ga+\rho)$ where $\rho\in\Lim$ and $F_i(\Ga,\rho)$ does not hold, which in turn is an uncritical case that recurs to equality of $\rho[m]$ and either $\eta[0]$
(where $\al$ is of a form $\thti(\De+\eta)$) or $\ual$. The situation where terms are of uncountable
cofinality is similar.
\end{rmk}

\begin{lem}\label{Ckvaluelem} Let $k>0$ and $\al\in\T\cap\Om_1$.
Then the sequence $(\Ck(\al[k-1]^i))_{i\le\istarkal}$ (with $\istarkal$ according to Definition \ref{iterateddescentdefi}) is weakly decreasing and passing through each $m\le\Ck(\al)$, 
such that for $i<\istarkal$ we have
\[\Ck(\al[k-1]^{i+1})<\Ck(\al[k-1]^i)\;\Longleftrightarrow\;\al[k-1]^i\mbox{  is a successor ordinal.}\] 
\end{lem}
{\bf Proof.} This is an immediate consequence of the definition of $\Ck$.
\qed

Recall the ordinals $\theta_n:=\thtnod(\ldots(\thtn(0))\ldots)$ from subsection \ref{prelimsec}, which form the ``master'' fundamental sequence of $\T\cap\Om_1$ that lies just outside
the system of fundamental sequences for $\T$ since $\T\cap\Om_1\not\in\T$, just as the ordinals $\om_n$ did (with respect to $\Esc$) in the classical Goodstein result for $\PA$.
 
\begin{cor}\label{maximalelementcor}
For any $k\ge2$ and $m<\om$ there exists a maximal $\al\in\T\cap\Om_1$ such that $\Ck(\al)=m$. Each $\theta_n$ is maximal in this sense, 
and the predecessor operation $\Pk$ preserves maximality (with respect to $k$).
\end{cor}
{\bf Proof.} The sequence $(\Ck(\tau_n))_{n<\om}$ is strictly increasing since $\tau_n=\tau_{n+1}[0]$, 
so that due to Bachmann property the sequence $(\tau_{n+1}[k-1]^i)_{i<\om}$ passes through $\tau_n$ after hitting successor stages, and moreover,
$\Ck(\be)>\Ck(\tau_n)$ for all $\be>\tau_n$. Thus each $\tau_n$ is maximal. 

Choosing $n$ such that $\Ck(\tau_n)>m$ and descending from $\tau_n$ by iteration of the operation $\cdot[k-1]$, according to Lemma \ref{Ckvaluelem} we eventually reach 
$\Ck(\tau_n[k-1]^i)=m$ at some step $i<\om$ for the first time, where again due to Bachmann property $\al:=\tau_n[k-1]^i$ is maximal such that $\Ck(\al)=m$. 
This also implies the claim regarding the predecessor operation $\Pk$.
\qed

\subsection{Quotients}\label{quotientssubsec}
We come to the central notion in this article, the notion of {\it quotient} of a set of ordinal terms. Corollary \ref{maximalelementcor} gives rise to the definition of a quotient 
$\TcapOmquok$ with respect to the equivalence relation induced by $\Ck$, i.e.\ $\al,\be\in\T\cap\Om_1$ are equivalent if and only if $\Ck(\al)=\Ck(\be)$, in which elements
can be identified with their maximal representative. The content of this subsection is elaborated to greater detail in \cite{W}.

\begin{defi}[Maximality quotients] \label{maximalityquotientdefi}For $k\ge2$ we define
\[\TcapOmquok:=\{\max\left\{\al\in\T\cap\Om_1\mid\Ck(\al)=m\}\mid m<\om\right\}.\]
\end{defi}

\begin{lem}\label{maxquotclslem} For $k\ge2$ the maximality quotient $\TcapOmquok$ is the closure of $\{\theta_n\mid n<\om\}$ under the predecessor operation $\Pk$. 
\end{lem}
{\bf Proof.} Immediate from 
Corollary \ref{maximalelementcor} and its proof. \qed

\begin{lem}\label{maxquotmonotolem} The sequence of maximality quotients $(\TcapOmquok)_{k\ge2}$ is $\subseteq$-increasing.
\end{lem}
{\bf Proof.} We have $\{\theta_n\mid n<\om\}\subseteq\TcapOmquok$ for every $k\ge2$. Any $\al\in\TcapOmquok$ is the result of iterated application of $\Pk$ to some $\theta_n$.
But, again as a consequence of Bachmann property, the sequence built by iterated applications of $\Pk$ to $\theta_n$ is a subsequence of the sequence built by iterations of $\Pke$ 
starting from $\theta_n$. So, we also have $\al\in\TcapOmquoke$.
\qed

\begin{theo}\label{maxquotisomtheo}
The restriction $\Ckstar$ of $\Ck$ to $\TcapOmquok$ \[\Ckstar:\quad\TcapOmquok\to\N\]
is an order isomorphism.
\end{theo}
{\bf Proof.} The sequences resulting from iterating $\Pk$ on the ordinals $\theta_n$ and writing them in increasing order results in end-extensions by finitely many elements
in each step of incrementation of $n$.
\qed

\begin{rmk}[cf.\ \cite{W}]\label{uniformityrmk} Note that maximality quotients can be abstractly established for notation systems $(\tau,\cdot[\cdot])$ where $\tau$ is a countable (recursive) ordinal and $\cdot[\cdot]$
is a system of fundamental sequences, for which we assume the following to hold:
\begin{enumerate}
 \item $0[n]=1[n]=0$ for any $n<\om$,
 \item $(\xi+\eta)[n]=\xi+\eta[n]$ for any $\xi<\tau$, $\eta\in(0,\tau)$, and $n<\om$,
 \item $(\la[n])_{n<\om}$, where $\la[0]>0$, is strictly increasing with supremum $\la$ for any limit ordinal $\la<\tau$, 
 \item \emph{Bachmann property:} For any $\al,\be<\tau$ and $n<\om$, if $\al[n]<\be<\al$, then $\al[n]\le\be[0]$, and
 \item there is a (strictly increasing) ``master'' fundamental sequence $(\theta_n)_{n<\om}$ such that $\sup\{\theta_n\mid n<\om\}=\tau$ 
 and $\theta_{n+1}[0]=\theta_n$ for all $n<\om$. 
 \end{enumerate}
The resulting order isomorphisms $\Ckstar$ for $\tau/_{\Ck}$ then give rise to a generalized Goodstein theorem of 
strength $\tau$ (as a proof-theoretic ordinal), with the generalization of Cichon's independence proof as described in the introduction and section \ref{Goodsteinprocesssec}, which reduces 
provability of the generalized Goodstein theorem to the provability of totality of $H_\tau$.  
\end{rmk}

\subsection{Quotients via restricted iteration}
In order to gain deeper insight into maximality quotients $\TcapOmquok$, we are going to characterize them by sets $\Tk$. To this end we start from the notion of a quotient $\Tquok$ 
of the entire set of notations $\T$, which will enable inductive proofs on the build-up of notations in $\T$. The equality of $\TcapOmquok$ and $\Tk$ will eventually follow from 
Lemma \ref{easyinclusionlem} and Corollary \ref{maincor}, see Theorem \ref{maintheo}.
 
\begin{defi}\label{quotientdefi}
We inductively define the elements of $(\Tquok)_{2\le k<\om}$ as increasing sequence of subsets of $\T$  as follows:
\begin{enumerate}
\item $0, 1\in\Tquok$.
\item If $\xi=_\ANF\xi_1+\ldots+\xi_n$ where $n\ge0$ and $\eta\in\Hz$ such that $\eta<\xi_n$ if $n>0$, then $\xi,\eta\in\Tquok$ implies that $\xi+\eta\cdot l\in\Tquok$
if and only if $l<k$.
\item If $\xi\in\T\cap\Omiz$, then $\al:=\thti(\xi)\in\Tquok$, provided that 
\begin{enumerate}
\item $\xipr\in\Tquok$ where $\xipr:=1+\xi$ if $i=0$ and $\xipr:=\xi$ otherwise,
\item $\al[0]\in\Tquok$, and 
\item $\al$ is not of a form $\be[n]$ where $\be\in\Tcirc\cap\Lim$ and $n\ge k-1$.
\end{enumerate}
\end{enumerate}
The quotients $\Tcircquok\subseteq\Tcirc$ are then defined by $\Tcircquok:=\Tquok\cap\Tcirc$, and we define the term sets relevant for the Goodstein process by
\[\Tk:=\Tcircquokbig\cap\Om_1.\]
\end{defi}

\begin{rmk}\label{clsrmk} $\mbox{ }$
\begin{enumerate}
\item Note the condition $1+j\in\Tquok$ for $\thtnod(j)=\om^{1+j}\in\Tquok$ where $j<\om$, which is why we introduced $\xipr$. 
\item Note that for any $\al\in\Tcircquok$ we expect $\al[j]\in\Tquok$ for $j<k-1$, so that the collapsed quotients $\Tk$ are closed under operations $\cdot[j]$ for $j<k-1$. 
This will eventually be verified by Lemma \ref{clslem}.
\item It is easy to see that we have \[\T=\bigcup_{k\in[2,\om)}\Tquokbig\quad\mbox{ and }\quad\Tcirc=\bigcup_{k\in[2,\om)}\Tcircquokbig.\]
\end{enumerate}
\end{rmk}

We are now going to show that the $<$-order type of the set 
$\Tk=\Tcircquok\cap\Om_1$ (where $k\in[2,\om)$) is $\om$ and that $\Ck$ enumerates this set, preserving its $<$-order.
To this end, we first define the analogue $\imc$ (for {\it iterative maximal coefficient}) of the measure $\mc$ that was used to define quotients $\Equok$ of $\Esc$.
Recall the $G$-norm from Definition \ref{iterateddescentdefi} and Corollary \ref{Gnormcor}. 

\begin{defi} The iterative maximal coefficient $\imc(\al)$ of terms $\al\in\T$ is defined recursively in $G(\al)$ as follows.
\begin{enumerate}
\item $\imc(0):=0$ and $\imc(1):=1$.
\item If $\xi=_\ANF\xi_1+\ldots+\xi_n\in\T$ where $n\ge0$ and $\eta\in\T\cap\Hz$ such that $\eta<\xi_n$ if $n>0$, then for $l\in(0,\om)$
\[\imc(\xi+\eta\cdot l):=\max\{\imc(\xi),\imc(\eta),l\}.\] 
\item If $\al\in\T$ is of a form $\thti(\xi)$, setting 
$\xipr:=\left\{\begin{array}{cl}
1+\xi & \mbox{  if }i=0\\
\xi     & \mbox{ otherwise}
\end{array}\right.$
we define
\[\imc(\al):=\max\{\imc(\xipr),\imc(\al[0]),n+1\},\] 
where either $n\ge1$ is maximal such that $\al$ is of a form $\al=\be[n]$ where $\be\in\Tcirc\cap\Lim$, or otherwise $n=0$.
\end{enumerate}
\end{defi}

\begin{rmk}
Note that we have $\imc(\al)=0$ if and only if $\al=0$. $\imc$ is weakly monotone with respect to the subterm relation, hence weakly increasing along 
localization sequences. 
\end{rmk}
 
We now obtain the following characterization of the quotients $\Tquok$ and $\Tcircquok$.
\begin{lem}
For $k\in[2,\om)$ we have \[\Tquok=\{\al\in\T\mid\imc(\al)<k\}\quad\mbox{ and }\quad\Tcircquok=\{\al\in\Tcirc\mid\imc(\al)<k\}.\]
\end{lem}
 {\bf Proof.} This follows immediately by induction along the involved inductive definitions.
\qed

\begin{lem}\label{intervalcollapsinglem} For $\la\in\Tcirc\cap\Lim$ and $\be\in\T\cap[\la[k-1],\la)$ where $k\ge2$ we have $\imc(\be)\ge k$ 
and hence \[\Tquok\cap[\la[k-1],\la)=\emptyset.\]  
\end{lem}
{\bf Proof.} We first observe that by definition of $\imc$, for any $\be\in\T$ \[\imc(\be)<k \Rightarrow \imc(\be[0])<k.\]
Suppose now that for some $\la\in\Tcirc\cap\Lim$, there exists a least $\be\in\T\cap[\la[k-1],\la)$ such that $\imc(\be)<k$. 
As $\imc(\la[k-1])\ge k$, this implies that $\la[k-1]<\be<\la$, with $\be\in\Lim$ since clearly $\imc(\ga)\le\imc(\ga+1)$ for any $\ga$. 
Bachmann property yields $\la[k-1]\le\be[0]$, which however contradicts the minimality of $\be$ since $\imc(\be[0])<k$. 
\qed

\begin{cor}[Maximality of elements in $\Tk$]\label{maximalitycor} 
Let $\al\in\Tk$. Then there does not exist any $\be\in\T\cap\Om_1$ such that $\Ck(\be)\le\Ck(\al)$ and $\be>\al$.
\end{cor}
{\bf Proof.} Suppose $\be\in\T\cap\Om_1$ satisfies $\Ck(\be)\le\Ck(\al)$. Recall Definition \ref{iterateddescentdefi}  and Lemma \ref{Ckvaluelem}, according to which the 
sequence $(\Ck(\be[k-1]^i))_{i\le\istarkbe}$ is weakly decreasing all the way down to $0$. Due to Lemma \ref{intervalcollapsinglem} $\al$ can not be element of any 
$[\be[k-1]^{i+1},\be[k-1]^i)$ for $i<\istarkbe$ such that $\be[k-1]^i\in\Lim$. Thus we must have $\al\ge\be$.
\qed

The following observations provide insight into how to verify that a term belongs to a quotient $\Tquok$.
\begin{lem}\label{fixpointinversionlem}
Let $\al=\thti(\De+\eta)\in\T$ be such that $\De+\eta\in\Tquok$, and denote its $\Omi$-localization by $(\Omi,\al_1,\ldots,\alm)$. 
Suppose that $\al=\ga[l+1]$ where $\ga\in\Tcirc\cap\Lim$ and $l\ge k-2$.
After ruling out the special case where $\al=\thtnod(k-1)$ and $\ga=\thtnod(\om)$, this has the following consequences:
\begin{enumerate}
\item $l=k-2$ and $\ga$ is of a form $\thti(\Ga+\rho)$ where $\chiomie(\Ga)=1$ and either $\rho\not\in\Lim$ or $F_i(\Ga,\rho)$ holds.
\item We either have $\eta=0$, or $\eta\in\Lim$ such that $F_i(\De,\eta)$ does not hold.
\item $\al$ can not be of a form $\be[0]$ where $\be=\thti(\Si+\nu)$, $\chiomie(\Si)=1$, and $F_i(\Si,\nu)$ holds.
\end{enumerate}
\end{lem}
{\bf Proof.} We verify the claims consecutively.
\\[2mm]
{\bf Ad 1.} Clearly, $\ga$ must be a $\thti$-term, and since $\De+\eta\in\Tquok$, Lemma \ref{inversionlem} shows that we can only have $\chiomie(\Ga)=1$ and either $\rho\not\in\Lim$
or $F_i(\Ga,\rho)$, as otherwise either $\De\not\in\Tquok$ or $\eta\not\in\Tquok$.  Moreover, only $l=k-2$ is possible.
\\[2mm]
{\bf Ad 2.} By part 1 we have $\De+\eta=\Ga[\ga[k-2]]$, and clearly $F_i(\Ga,\ga[k-2])$ does not hold.
\\[2mm]
{\bf Ad 3.} Since $\be[0]=\thti(\Si[\nu])$, assuming that $\al=\be[0]$, it follows that $\Ga[\ga[k-2]]=\De+\eta=\Si[\nu]$, hence either $\Si=\De+\Omie=\Ga$ and $\nu=\eta=\ga[k-2]$, or 
$\Ga,\De,\Si$ are limit-multiples of $\Omie$, $\Si[\nu]=\De=\Ga[\ga[k-2]]$, and $\eta=0$. Hence, using Corollary \ref{Bachmanncor}, $\De<\Ga=\Si$ and $\nu=\ga[k-2]$, but $\ga[k-2]$ 
has fixed point level below $\Ga=\Si$ despite $F_i(\Si,\nu)$: contradiction.
\qed

\begin{lem}\label{zeroiterationinversionlem}
Let $\al=\thti(\De+\eta)\in\T$ be such that $F_i(\De,\eta)$ holds, with $\Omi$-localization $(\Omi,\al_1,\ldots,\al_m)$, so that $\almmin=\eta$, $m>1$.
Suppose that $\De+\eta\in\Tquok$ and let $\hstar$ be such that $\al[0]^\hstar=\eta$.
Then for all $h\in(0,\hstar)$ the ordinal $\al[0]^h$ can not be of a form $\ga[l+1]$ where $\ga\in\Tcirc\cap\Lim$ and $l\ge k-2$.
\end{lem}
{\bf Proof.} Note that $\hstar$ exists due to part 2 of Lemma \ref{supportbyzeroiterationlem}. As the case $\De=0$ is trivial ($\hstar=1$), parts 4 (including its proof) and 5 of 
Lemma \ref{supportbyzeroiterationlem} together with Lemma \ref{fixpointinversionlem} apply to $\al[0]^h$, $h\in(0,\hstar)$, and yield the claim. 
\qed

\section{Characterization theorem for quotients}
The upcoming Theorem \ref{maintheo} is in correspondence with  Lemma 1 and Remarks 1 and 2 of \cite{Cichon1983}. Since the generalized Goodstein principle introduced here induces
much faster growing operations on the natural numbers than just base-$k$ exponentiation, Lemma 1 of \cite{Cichon1983} does not generalize directly. It is therefore instructive
to sketch an alternative proof for the fact that the mapping $\Ck$ restricted to the quotient $\Equok$ is the Mostowski collapse of $\Equok$ onto $\om$.
This proof, which follows from the next lemma, then does generalize to the quotients $\Tk$. 
For limit $\al$, by $i_\al$  we count how many times we can iterate the operation $\cdot[k-1]$ while staying within $\Lim$.

\begin{defi}\label{alternativedefi} 
For $k\in[2,\om)$ and $\al\in\Esc\cap\Lim$ let 
\begin{equation}\label{ialeq}
i_\al:=\max\left\{i<\om\mid\forall j\le i\:\left(\al[k-1]^j\in\Lim\right)\right\}.
\end{equation}
\end{defi}
The ordinal $\al[k-1]^{i_\al}$ then clearly is a successor-multiple of $\om$, so that $\al_0:=\al[k-1]^{i_\al}-\om$ is well-defined and \[\Ck(\al)=\Ck(\al[k-1]^{i_\al+1})=\Ck(\al_0)+k.\]
\begin{lem}\label{alternativelem} 
For $k\in[2,\om)$ and $\al\in\Equok\cap\Lim$ we have \[\al_0\in\Equok.\]
\end{lem}
{\bf Proof.} This now follows by $<$-induction on $\al\in\Equok$. If $\al=_\NF\xi+\eta$, the claim follows directly from the i.h.\ for $\eta$. If $\al=\om^{\be+1}$ for some $\be$,
we have $\al_0=0$ in case of $\be=0$. Otherwise we have $\al[k-1]=\om^\be\cdot(k-1)+\om^\be$, apply the i.h.\ to $\om^\be$, and the claim follows.
If finally $\al=\om^\la$ for some $\la\in\Lim$, by the i.h.\ for $\la$ we obtain $\la_0:=\la[k-1]^{i_\la}-\om\in\Equok$. In case of $\la_0=0$, we must have $\la=\om$ and hence $i_\la=0$.
It follows that $i_\al=k$ and \[\al_0=(\om^\om)_0=\om^{k-1}\cdot(k-1)+\ldots+\om\cdot(k-1)\in\Equok.\] Otherwise we have \[\al[k-1]^{i_\la+k+1}=\om^{\la_0+k-1}\cdot(k-1)+\ldots+\om^{\la_0}\cdot(k-1)+\om^{\la_0}\]
and obtain the claim with $i_\al=i_\la+k+1+i_{\om^{\la_0}}$ after an application of the i.h.\ to $\om^{\la_0}$.
\qed

\begin{cor}\label{esccollapscor} For $k\ge2$ the mapping $\Ck$ restricted to the quotient $\Equok$ is the Mostowski collapse of $\Equok$ onto $\N$
and the predecessor function $\Pk$ restricted to $\Equok$ is the true predecessor function on $\Equok$.
\end{cor}
{\bf Proof.} As in Lemma \ref{intervalcollapsinglem} for quotients $\Tquok$, we have a corresponding situation in $\Equok$: For any $\la\in\Esc\cap\Lim$ and $\be\in[\la[k-1],\la)$
we have $\mc(\be)\ge k$, thus \[\Equok\cap[\la[k-1],\la)=\emptyset.\]
By Lemma \ref{alternativelem}, for any $\al\in\Equok\cap\Lim$ the mapping $\Ck$ counts the elements of 
\[[\al_0,\al]\cap\Equok=\{\al_0\}\cup\{\al[k-1]^{i_\al}[j]\mid j<k-1\}\cup\{\al\}\] starting from $\Ck(\al_0)$ to $\Ck(\al)$ in increasing order.
Note that for $\al\in\Equok\cap\Lim$ we have $\Pk\al=\al_0+k-1$.
 \qed

We now transfer this argumentation to terms in $\T$. Recall that $\sumend(\ga)$ of an ordinal $\ga$ denotes the last additive component of $\ga$.  
Clearly, as seen by straightforward $<$-induction on $\al$, iterated application of the operation $\cdot[j]$, $j<\om$, to any $\al\in\T\cap\Hz^{>1}$ eventually reaches a 
successor-multiple of $\Oml$ where $l:=\lv(\al)$. However, for terms $\al\in\Tcirc$ such that $\lv(\al)>0$, iterated application of $\cdot[j]$, $j<\om$, is only meaningful until
a term $\be$ such that $\domf(\be)>0$ is reached for the first time. 
\begin{defi}\label{descentdefi} Let $k\ge2$, $\al\in\T\cap\Lim$, $l:=\lv(\sumend(\al))$, and define $i_\al$, extending (\ref{ialeq}), by
\begin{equation}\label{ialexteq}
i_\al:=\left\{\begin{array}{cl}
\max\left\{i<\om\mid\forall j\le i\:\left(\al[k-1]^j\in\Lim\right)\right\} & \mbox{ if } l=0\\[2mm]
\min\left\{i<\om\mid\domf(\al[k-1]^i)>0\right\} & \mbox{ if } l>0,
\end{array}\right.
\end{equation}
and $\al_0$ by 
\[\al_0:=\left\{\begin{array}{cl}
\al[k-1]^{i_\al}-\om & \mbox{ if } l=0\\[2mm]
\al[k-1]^{i_\al} & \mbox{ if } l>0.
\end{array}\right.\]
\end{defi}

\begin{lem}\label{easyinclusionlem} Let $k\ge2$. We have \[\Tk\subseteq\TcapOmquok.\]
\end{lem}
{\bf Proof.} Recall that $\TcapOmquok$ is the closure of the set $\{\theta_n\mid n<\om\}$ under the predecessor operation $\Pk$, see Corollary \ref{maxquotclslem}.
Clearly, $\{\theta_n\mid n<\om\}\subseteq\Tk$ as $\imc(\theta_n)=1$ for each $n<\om$, and $\sup\{\theta_n\mid n<\om\}=\sup\Tk=\T\cap\Om_1$. 
For each $\al\in\Tk\cap\Lim$, by Lemma \ref{intervalcollapsinglem} we have \[[(\al)_0+k,\al)\cap\Tk=\emptyset.\]
Moreover, for $\al\in\Tk\cap\Lim$ the ordinal $\al_0$ enjoys the maximality property that according to Corollary \ref{maximalitycor} elements of $\Tk$ do:
Assume the existence of some $\be\in\T\cap\Om_1$ such that $\be>\al_0$ and $\Ck(\be)\le\Ck(\al_0)$. 
As a consequence of the definition of $\Ck$ we have $\Ck(\al)=\Ck(\al_0)+k$, and due to Bachmann property, for all $\ga\in(\al_0,\al]$ we have \[\Ck(\ga)>\Ck(\al_0).\]
It follows that $\be>\al$, which however contradicts Corollary \ref{maximalitycor} since $\Ck(\be)\le\Ck(\al_0)<\Ck(\al)$. 
\qed

The reverse inclusion $\TcapOmquok\subseteq\Tk$ will follow from the Main Lemma of this article, Lemma \ref{mainlem}, which establishes that $\al_0\in\Tk$ for limits $\al\in\Tk$ and hence
that $\Tk$ is closed under $\Pk$, see Corollary \ref{maincor}.
Once Bachmann property for the system of fundamental sequences is established as in \cite{W26}, this turns out to be the most involved part of the argumentation.

Since the computation of $\Ck(\al)$ for $\al\in\Tk$ involves intermediate calculation of $\Ck(\be)$ for some $\be\in\Tcirc\setminus\Tk$, we extend the set $\Tcircquok$ minimally so as 
to include these auxiliary ordinals (terms) to enable an inductive proof. In fact, we need such extension only for terms in $\Tcircquok\cap\Hz$. 
This extension can be seen as a suitable closure of the set $\Tcircquok\cap\Hz$.
Recall that $\Om_0^\prime:=1$ and $\Omipr=\Omi$ for any $i>0$.

\begin{defi}\label{Tcircquokclsdefi}
Let $k\ge2$. We define the following sets of terms:
\begin{enumerate}
\item $\al\in\T$ is called {\it $k$-residual} if $\al=\ga[k-1]^j$ for some $\ga\in\Tcircquok$ and $j\le i_\ga$. 
\item $\Tscircquok:=\{\al=\thti(\De+\etapr+1)\mid\al\mbox{ is $k$-residual and }\ual\in\Tcircquok\}$.
\item $\Tfcircquok:=\{\al=\thti(\De+\eta)\mid(\eta\not\in\Lim\vee F_i(\De,\eta))\wedge\chiomie(\De)=1 \mbox{ such that } \al \mbox{ is k-tangent}\}$,
where such fixed-point terms $\al$ are called {\it $k$-tangent} if the following conditions apply:
\begin{enumerate}
\item $\al\in\Tcirc$ is $k$-residual,
\item $\ual\in\Tquok$,
\item $\De[\be]\in\Tquok$ for any $\be\in\Tquok\cap\Hz\cap[\Omipr,\Omie)\cup\{0\}$, 
\item $\al[0],\ldots,\al[k-2]\in\Tquok$.
\end{enumerate}
\item $\Tlcircquok:=\{\ga[k-1]\mid\ga\in\Tfcircquok\mbox{ and }\ga[k-1]\in\Tcirc\}$.
\item The closure $\Tcircquokcls$ of $\Tcircquok\cap\Hz$ by additive principal terms occurring during calculation of $\Ck\!\restriction_\Tk$ is now defined by
\[\Tcircquokcls:=(\Tcircquok\cap\Hz)\cup\Tscircquok\cup\Tfcircquok\cup\Tlcircquok.\]
\end{enumerate}
\end{defi}

\begin{lem}\label{supporttermlem} We have $\Tcircquokcls\subseteq\Tcirc\cap\Hz$, $\alcirc, \ual\in\Tquok$ for all $\al\in\Tcircquokcls$, and $(\al[k-1])^\circ\in\Tquok$ if $\al\in\Tfcircquok$
also if $\al[k-1]\not\in\Tcirc$.
\end{lem}
{\bf Proof.} By part 1 of Lemma \ref{supportbyzeroiterationlem} and Corollary \ref{Gnormcor} the quotient $\Tquok$ is closed under the operations $\al\mapsto\ual$ and $\al\mapsto\alcirc$.
If $\al\in\Tcircquok\cup\Tscircquok$ we are done. Assume from now that $\al\in\Tfcircquok$.
We then observe that for $\al[k-1]$ we have $\al[k-1]\in\Tlcircquok$ only if $\al[k-1]\in\Tcirc$, which, by the way, excludes
terms of the form $\thti(\Omie)[1]=\Om_i^2=(\thti(\Omie))_0$ where $i>0$, and we also have $\underline{\al[k-1]}=0$:
Indeed, for $\al$ of the form $\thti(\De+\eta)$ such that $\chiomie(\De)=1$ we have $(\De[\al[k-2]])^\stari=\al[k-2]$, which is not a successor nor
can it be equal to $\ual$ or to the predecessor of $\al[k-1]$ in its $\Omi$-localization. In case of $\De=\Depr+\Omie$ for some $\Depr$ we have $\neg F_i(\Depr,\al[k-2])$. 

We finally show that $\alcirc, (\al[k-1])^\circ\in\Tquok$ for $\al=\thti(\De+\eta)\in\Tfcircquok$. Technically, $6$ combined cases occur: writing $\De=_\NF\Depr+\sumend(\De)$, we either 
have (A) $\sumend(\De)=\Omie$ or (B) $\sumend(\De)$ is a limit-multiple of $\Omie$ such that $\chiomie(\sumend(\De))=1$, and for $\eta$ we either have (a) $\eta=\etapr+1$,
(b) $F_i(\De,\eta)$ holds, or (c) $\eta=0$. In situation (a) we have $\alcirc=\ual=\thti(\De+\etapr)\in\Tquok$, in situation (A) we have $(\al[k-1])^\circ=\al[0]\in\Tquok$, which in situations
(Ab) and (Ac) is equal to $\alcirc$, and in situation (B) we have $(\al[k-1])^\circ=(\al[0])^\circ\in\Tquok$, which in situations (Bb) and (Bc) is equal to $\alcirc$. Note that in situation (c)
we distinguish whether $m=1\vee\almmin\le(\Depr)^\stari$, resulting in $\alcirc=\thti(\Depr)$, or otherwise, resulting in $\alcirc=\thti(\Depr+\almmin)$. 
\qed

\begin{lem}\label{Tquoksubstlem}
Suppose $\Ga=\thtje(\Xi+\rho)\in\Tquok$ satisfies $\chiomie(\Ga)=1$ for some $i$. Then we have \[\Ga[\be]\in\Tquok\]
for any $\be\in\Tquok\cap\Hz\cap[\Omipr,\Omie)\cup\{0\}$.
\end{lem}
{\bf Proof.} The proof is by induction on $G(\Ga)$ if $i>0$, and on $G(\Ga[1])$ if $i=0$. Recall Corollary \ref{Gnormcor}. 
Note that the assumptions imply $i\le j$ and
that $\Ga[0]\in\Tquok$ is immediate by assumption and definition of $\Tquok$. We may therefore assume that $\be>0$. 
We will treat the situation $i=0$ and $\be=1$ as special but crucial case. The regular case is where $\be\in\Hz^{>1}$. 
Let $(\Omje,\Ga_1,\ldots,\Ga_m)$ be the $\Omje$-localization of $\Ga$. 
\\[2mm]
{\bf Case 1:} $\rho\in\Lim$ and $F_{j+1}(\Xi,\rho)$ does not hold. Then $\chiomie(\Ga)=\chiomie(\rho)=1$, $\Ga[\be]=\thtje(\Xi+\rho[\be])$, and by the i.h.\ for $\sumend(\rho)$ we have 
$\rho[\be]\in\Tquok$, hence also $\Xi+\rho[\be]\in\Tquok$.
According to Lemma \ref{zeroiterationlem} we have \[\Ga[\be][0]^h=\thtje(\Xi+\rho[\be[0]^h])=\Ga[\be[0]^h]\] for $h$ such that $\be[0]^{h^\prime}\in\Lim$ for all $h^\prime<h$, 
and by the i.h.\ for $\rho$ we have $\rho[\be[0]^h]\in\Tquok$. Let $\hstar$ be the maximal such $h$. We then have $\be[0]^\hstar=1$ if $i=0$ and $\be[0]^\hstar=0$ otherwise, thus
\[\Ga[\be][0]^\hstar=\Ga[\be[0]^\hstar]=\left\{\begin{array}{cl} 
\thtje(\Xi+\rho[1]) = \Ga[1] & \mbox{ if } i=0\\[2mm]
\thtje(\Xi+\rho[0]) = \Ga[0] & \mbox{ if } i>0.
\end{array}\right.\] 
Note that for $\be\in[\Omipr,\Omie)$, according to Corollary \ref{Bachmanncor} $\Ga[\be]$ can not be of a form $\ga[l+1]$ for any $\ga\in\Tcirc\cap\Lim$, $l<\om$. 
Thus, in case of $i>0$ we are done, as $\Ga[0]\in\Tquok$. Otherwise, we are left with showing that $\Ga[1]\in\Tquok$ which in turn reduces to verifying that $\Ga[1][0]\in\Tquok$.
\\[2mm]
{\bf Subcase 1.1:} $\rho=_\NF\rhopr+\Om_1$. Then we have $\rho[1]=\rhopr+1$, so that $\Ga[1]=\thtje(\Xi+\rhopr+1)$,
$\Gapr:=\underline{\Ga[1]}=\thtje(\Xi+\rhopr)=\Ga[0]\in\Tquok$, and
\[\Ga[1][0]=\left\{\begin{array}{cl}
\Gapr & \mbox{ if }\Xi=0\\[2mm]
\thtje(\Xi[\Gapr]) & \mbox{ if }\chiomjz(\Xi)=1\\[2mm]
\thtje(\Xi[0]+\Gapr) & \mbox{ otherwise.}
\end{array}\right.\]
We are done if $\Xi=0$, so assume that $\Xi>0$ from now. If $\chiomjz(\Xi)=0$, we have $\Xi[0]+\Gapr\in\Tquok$, as $\Xi[0]$ is still a multiple of $\Omjz$, whereas if
$\chiomjz(\Xi)=1$, the i.h.\ applies to $\Xi$, noting that $G(\Xi)<G(\Ga[1])$ and $\Gapr\ge\Om_1$, hence $\Xi[\Gapr]\in\Tquok$.
As $\Gapr$ is a non-zero support term, Lemmas \ref{fixpointinversionlem} and \ref{zeroiterationinversionlem} show that terms in $\zeroseq(\Ga[1][0])$ above $\Gapr$
can not be of a form $\ga[l+1]$ for any $\ga\in\Tcirc\cap\Lim$, $l\ge k-2$.
\\[2mm]
{\bf Subcase 1.2:} $\rho$ is a limit-multiple of $\Om_1$. Then we have 
\[\Ga[1][0]=\thtje(\Xi+\rho[1][0]),\]
where $\rho[1][0]\in\Tquok$ since $\rho[1]\in\Tquok$ and $\rho[1]$ is a limit-multiple of $\Om_1$ according to Lemma \ref{zeroiterationlem}. 
Recalling Definition \ref{socletermdefi}, we see that $\Ga_{m-1}\le\Ga^\circ\in\zeroseq(\Ga)$ where 
\[\Ga^\circ=\left\{\begin{array}{cl}
\thtje(\Xi+\rhopr) & \mbox{ if }\rho=_\NF\rhopr+\sumend(\rho)\not\in\Hz,\\[2mm]
\thtje(\Xi+\Ga_{m-1}) & \mbox{ if }\Xi^\starje<\Ga_{m-1}<\rho \mbox{ with }m>1,\\[2mm] 
\thtje(\Xi) & \mbox{ otherwise.}
\end{array}\right.\]
We observe that $\Ga^\circ\in\zeroseq(\Ga[1][0])$, and are therefore left with showing that none of the terms in $\zeroseq(\Ga[1][0])$ above $\Ga^\circ$, which consist of all terms 
\[\Ga[1][0]^h=\thtje(\Xi+\rho[1][0]^h)\] for $h\in(0,\hstar]$ where $\hstar$ is the maximal $h$ such that $\Ga[1][0]^h>\Ga^\circ$,
can be of the form $\ga[l+1]$ where $\ga\in\Tcirc\cap\Lim$ and $l\ge k-2$. 
Since $\Xi, \rho[1][0]^h\in\Tquok$ and $\rho[1][0]^h>0$ for $0<h\le\hstar$, this would imply that $\ga$ is of a form $\thtje(\Xi+\Omje+\si)$ for a suitable $\si$ and that $l=k-2$, 
so that $\rho[1][0]^h=\ga[k-2]$.
But then we would have $\ga[k-1]=\Ga[1][0]^h<\Ga<\ga$, as $\ga\in(\Ga_{m-1},\Ga)$ is impossible since $\ga$ has a strictly greater fixed-point level than $\Ga$, see Proposition \ref{localipic},
and Lemma \ref{intervalcollapsinglem} would entail $\Ga\not\in\Tquok$, contradiction. 
\\[2mm]
{\bf Case 2:} $\rho\not\in\Lim$ or $F_{j+1}(\Xi,\rho)$ holds. Since $\chiomie(\Ga)=\chiomie(\Xi)=1$ and $i\le j$, we have $\Xi>0$ and $\chiomjz(\Xi)=0$.
Thus $\Ga[\be]=\thtje(\Xi[\be]+\uGa)$, where $\uGa\in\Tquok$ and $\Xi[\be]\in\Tquok$ by the i.h.\ for $\sumend(\Xi)$.
Again according to Lemma \ref{zeroiterationlem} we have \[\Ga[\be][0]^h=\thtje(\Xi[\be[0]^h]+\uGa)=\Ga[\be[0]^h]\] for $h$ such that $\be[0]^{h^\prime}\in\Lim$ for all $h^\prime<h$, 
and by the i.h.\ for $\sumend(\Xi)$ we have $\Xi[\be[0]^h]\in\Tquok$. Let $\hstar$ be the maximal such $h$. We thus have 
\[\Ga[\be][0]^\hstar=\Ga[\be[0]^\hstar]=\left\{\begin{array}{cl} 
\thtje(\Xi[1]+\uGa) = \Ga[1] & \mbox{ if } i=0\\[2mm]
\thtje(\Xi[0]+\uGa) = \Ga[0] & \mbox{ if } i>0.
\end{array}\right.\] 
As in Case 1, for $\be\in[\Omipr,\Omie)$ Corollary \ref{Bachmanncor} implies that $\Ga[\be]$ can not be of a form $\ga[l+1]$ for any $\ga\in\Tcirc\cap\Lim$, $l<\om$.
If $i>0$, we are done since $\Ga[0]\in\Tquok$. 

In case of $i=0$, we have to show that $\Ga[1]\in\Tquok$. According to Lemma \ref{zeroiterationlem} $\Xi[1]$ is a limit-multiple of $\Omjz$ and $\domf(\Xi[1])\le j+2$, 
and by Lemma \ref{localizationlem} the immediate predecessor of $\Ga[1]$ in its $\Omje$-localization is $\uGa=\thtje(\Xi+\rhopr)$ in case of $\rho=\rhopr+1$ and $\Ga_{m-1}$
otherwise. We have 
\[\Gapr:=\underline{\Ga[1]}=\left\{\begin{array}{cl}
\uGa & \mbox{ if }\uGa>0, \mbox{ otherwise:}\\[2mm]
\Ga_{m-1} & \mbox{ if }\Xi[1][0]^\starje<\Xi[1]^\starje=\Ga_{m-1}\mbox{ where }m>1,\\[2mm]
0 & \mbox{ otherwise,}
\end{array}\right.\] 
hence $\Gapr\in\Tquok$, and note that $\Gapr=\uGa$ if $\chiomjz(\Xi[1])=1$.
We thus have 
\[\Ga[1][0]=\left\{\begin{array}{cl}
\thtje(\Xi[1][\Gapr]) & \mbox{ if } \chiomjz(\Xi[1])=1\\[2mm]
\thtje(\Xi[1][0]+\Gapr) & \mbox{ otherwise,}
\end{array}\right.\]
and note that if $\chiomjz(\Xi[1])=1$, the i.h.\ applies to $\Xi[1]$ since $G(\Xi[1])<G(\Ga[1])$ and $j+1>0$ (noting that $j+1$ takes the role of $i$ in this instance), showing that $\Xi[1][\Gapr]\in\Tquok$. 

Now, if $\Gapr>0$, we argue as before using Lemmas \ref{fixpointinversionlem} and \ref{zeroiterationinversionlem} to see that terms in $\zeroseq(\Ga[1][0])$ above $\Gapr$ can not be of the form 
$\ga[l+1]$. Assume, therefore, that $\Gapr=0$ from now on. We then have $\Ga[1][0]=\thtje(\Xi[1][0])$ and need to show that terms $\Ga[1][0]^h$ for $h\in(0,\hstar]$, where $\hstar$ is the maximal
$h$ such that $\Ga[1][0]^h>\Ga_{m-1}$, can not be of the form $\ga[l+1]$. In case there exists a minimal $\hpr\le\hstar$ such that $(\Xi[1][0]^\hpr)^\starje<\Ga_{m-1}$, we argue as above 
(with $\Ga_{m-1}$ in place of $\Gapr$) for $h\ge\hpr$, so we are left with checking terms $\Ga[1][0]^h=\thtje(\Xi[1][0]^h)$ for $h\in(0,\hpr)$.

Assuming that $\ga[l+1]=\Ga[1][0]^h=\thtje(\Xi[1][0]^h)\in(\Ga_{m-1},\Ga)$ for some $\ga\in\Tcirc\cap\Lim$ and $l\ge k-2$ implies that $\ga$ is of a form $\ga=\thtje(\Si+\si)$ and $l=k-2$, 
since clearly $\Xi[1][0]^h\in\Tquok$. We thus have $\Xi[1][0]^h=\Si[\ga[k-2]]$. Bachmann property implies that $\Xi<\Si$ (using Corollary \ref{Bachmanncor} to see first that 
$\Xi[1][0]^{h-1}<\Si$), so that by the property of localization, see Proposition \ref{localipic}, 
$\ga\not\in(\Ga_{m-1},\Ga)$, hence $\Ga<\ga$. But since $\ga[k-1]=\Ga[1][0]^h<\Ga<\ga$, Lemma \ref{intervalcollapsinglem} then yields the contradiction 
$\Ga\not\in\Tquok$.
\qed

\begin{lem}\label{Tquoksimplefixplem}
Let $\al=\thti(\De+\eta)\in\Tcircquok\cup\Tscircquok$ be such that $\chiomie(\De)=1$ and either $\eta\not\in\Lim$ or $F_i(\De,\eta)$ holds. 
Then we have
\[\al[0],\ldots,\al[k-2]\in\Tquok.\]
\end{lem}
{\bf Proof.} The assumptions imply that $\De,\ual\in\Tquok$, and if $\al\in\Tcircquok$, we have $\al[0]\in\Tquok$ by definition of $\Tquok$. 
Otherwise we have $\ual=\thti(\De+\etapr)\in\Tcircquok$ where $\eta=\etapr+1$.
According to Lemma \ref{Tquoksubstlem} we have $\De[\be]\in\Tquok$ for all $\be\in\Tquok\cap\Hz\cap[\Omipr,\Omie)\cup\{0\}$.
By induction on $j$ we now verify that $\al[j]\in\Tquok$ for $j\le k-2$. Let $\alpr:=\ual$ if $j=0$ and $\alpr:=\al[j-1]$ otherwise. We then have $\De[\alpr]\in\Tquok$ using the i.h.\ for $j>0$.
We need to show that $\al[j][0]^h\in\Tquok$ for any $h$, which includes showing that $\al[j][0]^h$ can not
be of the form $\ga[l+1]$ where $\ga\in\Tcirc\cap\Lim$ and $l\ge k-2$. For $h=0$ this directly follows from Corollary \ref{Bachmanncor} if $j>0$ and is immediate for $j=0$ if $\al\in\Tquok$.

Otherwise, if $j=0$ and $\al\not\in\Tquok$, we first show that $\al[0]\in\Tquok$. Assuming that $\thti(\De[\ual])=\ga[l+1]=\thti(\Ga+\rho)[l+1]$, since $\De[\ual]\in\Tquok$, it follows that 
$\chiomie(\Ga)=1$, hence $\De[\ual]=\Ga[\ga[l]]$, which 
is impossible since $\ual,\ga[l]>0$, noting that $\Ga=\De$ does not occur as $\ual$ is of greater fixed-point level than $\ga[l]$, and applying Corollary \ref{Bachmanncor} otherwise. 
Iterating the application of $\cdot[0]$ to $\al[0]=\thti(\De[\ual])$ eventually reaches $\ual\in\Tquok$ itself, where in case of $i=0$ we pass through the term $\thtnod(\De[1]+\ual)$. 
Using Lemmas \ref{zeroiterationlem} and \ref{zeroiterationinversionlem}, or directly by side induction on $G(\al[0]^h)-G(\ual)$ one readily shows that none of the terms 
$\thti(\De[\ual])[0]^{h+1}$ can be of a form $\ga[l+1]$ where $\ga\in\Tcirc\cap\Lim$ and $l\ge k-2$, as the support term $\ual$, which has fixed-point level $\De$, 
remains constant throughout and would have to match $\ga[k-2]$, which however is of a strictly lower fixed-point level.  
Thus, $\al[0]\in\Tquok$ also if $\al\in\Tscircquok\setminus\Tcircquok$.

{\bf Inductive step:}
Suppose that $\al[j]\in\Tquok$ for some $j<k-2$. We have to show that $\al[j+1]=\thti(\De[\al[j]])\in\Tquok$. Since $\al[j]\in\zeroseq(\al[j+1])$ we only need to consider terms 
$\al[j+1][0]^h$ above $\al[j]$ and argue by side induction on $G(\al[j+1][0]^h)-G(\al[j])$.
Indeed, using Lemmas \ref{zeroiterationlem} and \ref{supportbyzeroiterationlem} \[\al[j+1][0]^h=\thti(\De[\al[j][0]^h])=\ldots=\thti(\De[\ldots\thti(\De[\al[0]^{h+1}])\ldots]),\] 
where $h$ is such that $\al[0]^{h+1}\ge\ual$. More explicitly, if $\ual>0$, we have $\al[0]^{h+2}=\thti(\De[\ual[0]^{h+1}]+\ual)$, which in case of $i=0$ intermediately passes through the
term $\thti(\De[1]+\ual)$, and eventually reaches $\ual\in\Tquok$. If $\ual=0$, let $(\Omi,\al_1,\ldots,\alm)$ be the $\Omi$-localization of $\al$ and note that if $m>1$,
$\almmin\le\De[0]^\stari$ since $\chiomie(\De)=1$. Terms $\al[0]^{h+1}$ then eventually pass through
$\almmin,\ldots,\al_1$ if $m>1$, using these terms as support terms before reaching them, and when descending from $\al_1$ eventually reach $\Omi$.
Note that in the situation $i=0=\ual$ intermediately the subterm $\al[1][0]^\hpr=\thtnod(\De[1])$ occurs, from which further iteration of $\cdot[0]$
leads to $\al[0]=\thtnod(\De[0])$. Thus there exists $\hstar$ such that \[\al[j+1][0]^\hstar=\al[j].\] 
Now, if we attempt to match $\al[j+1][0]^h$, where $h<\hstar$, with a term $\ga[l+1]$, we again quickly find that
this would first imply that $\ga$ is of a form $\thti(\Ga+\rho)$ with $\chiomie(\Ga)=1$ and suitable $\rho$, as well as that $l=k-2$, since $\al[j+1][0]^h$ is a $\thti$-term the argument of 
which is an element of $\Tquok$. Using Corollary \ref{Bachmanncor}, we observe next that this would further imply that $\Ga=\De$, and since $j<k-2$ parsing terms leads to a 
cancellation of the leading $j$-many nestings of $\thti(\De[\cdot])$ until the equality of \[\al[1][0]^h=\ga[k-j-2]=\thti(\De[\ga[k-j-3]])\] is obtained. 
We have $\al[1][0]^h=\thti(\De[\al[0]^{h+1}])$, and if $\ual>0$, we argue similarly as in the verification of $\al[0]\in\Tquok$ above. Suppose therefore that $\ual=0$ from now. 
Iterating $\cdot[0]$ on $\al$ leads to $1$ in case of $i=0$, otherwise to $0$, hence $\al[1][0]^h$ eventually reaches $\thti(\De[1])$ if $i=0$ and directly $\al[0]$ otherwise.

If $i>0$, intermediate terms $\al[0]^{h+1}$ are of a form either $\thti(\De[0]^r)$, $\thti(\De[0]^r+\al_s)$, or $\thti(\De[0]^r[\al_s])$, 
along with $\thti(\De[0]^r[\al_s[0]^p]+\al_s)$ and eventually $\al_s$, where $r\ge2$ and $s\in\{1,\ldots,m-1\}$. In the cases where an $\al_s$ occurs explicitly, matching with $\ga$ fails immediately.
This leaves us with the situation $\thti(\De[0]^r)=\ga[k-j-3]=\thti(\De[\gapr])$ where $\gapr$ is either $\uga$ or $\ga[k-j-4]$. Since $r\ge2$ this is impossible as well.

If $i=0$, we argue as before for terms $\al[1][0]^h$ until we reach the term $\al[1][0]^\hstar=\thti(\De[1])$. But then we need to have $\thti(\De[1])=\thti(\De[\ga[k-j-3]])$ where $k-j-3\ge0$,
which is impossible. Clearly, descending further via $\cdot[0]$ from $\thti(\De[1])$ does not match with $\thti(\De[\ga[k-j-3])$ either.
\qed

Although it is not directly required to establish Lemma \ref{mainlem} below, we state the expected closure property mentioned in part 2 of Remark \ref{clsrmk}.

\begin{lem}\label{clslem} 
For any $\al\in\Tcircquok$ we have $\al[0],\ldots,\al[k-2]\in\Tquok$.
\end{lem}
{\bf Proof.} Clearly, we may assume that $\al$ is of a form $\thti(\De+\eta)$ and proceed by induction on the build-up of $\al$ and side induction on $j\le k-2$, where the base case 
$j=0$ comes for free by the definition of $\Tquok$.

The special cases $\al\le\om^\om$ are seen immediately, and the usual case distinction applies then for $\al>\om^\om$. 
If $\eta\not\in\Lim$ or $F_i(\De,\eta)$ holds, the case $\De=0$ follows directly from $\ual\in\Tquok$ according to Lemma \ref{supporttermlem}, 
Lemma \ref{Tquoksimplefixplem} handles the situation where $\chiomie(\De)=1$,
and if $\De>0$ and $\chiomie(\De)=0$, the ``$\De$-limit case'', we have $\domf(\De)=\domf(\al)=0$, so that the i.h.\ applies to $\De$ and $\al[j+1]=\thti(\De[j+1]+\ual)$.
If in this situation $\ual>0$, Lemmas \ref{fixpointinversionlem} and \ref{zeroiterationinversionlem} apply, while for $\ual=0$ the usual localization argument applies to the elements
of $\zeroseq(\thti(\De[j+1]))$ above $\almmin$ where $(\Omi,\al_1,\ldots,\alm)$ is the $\Omi$-localization of $\al$.

The ``$\eta$-limit case'', that is, $\eta\in\Lim$ such that $F_i(\De,\eta)$ does not hold, finally applies the i.h.\ to $\eta$ and the side i.h.\ to $\al[j]$ and considers $\zeroseq(\al[j+1])$ 
above the eventually reached $\al[j]$, with special treatment of the situation $\sumend(\eta)=\om$, where Lemma \ref{Tquoksubstlem} applies to see that $\De[\al[j]]\in\Tquok$ 
in case of $\chiomie(\De)=1$.
\qed

Recall the operation $\al\mapsto\al_0$ from Definition \ref{descentdefi}, which will be used frequently in the upcoming lemmas.

\begin{lem}\label{Tquoklimlem}
Suppose that $\al=\thti(\De+\eta)$, where $\eta\in\Lim$ and $F_i(\De,\eta)$ does not hold, is $k$-residual such that $\alcirc, \De+\eta\in\Tquok$. Let $l:=\lv(\sumend(\eta))$ and $\eta_0$ be
according to Definition \ref{descentdefi}.
\begin{enumerate}
\item If $l=0$, so that $\eta[k-1]^{i_\eta}=_\NF\eta_0+\om$, and $\eta_0\in\Tcircquok$, then setting $\alpr:=\al[k-1]^{i_\eta+1}$ we have 
\[\ualpr\in\Tcircquok\]
and therefore $\alpr\in\Tscircquok$.
\item If $l>0$, so that $\eta[k-1]^{i_\eta}=\eta_0$ where $\domf(\eta_0)=j+1$ for some $j$, and provided that $\eta_0[\be]\in\Tquok$ for all $\be\in\Tquok\cap\Hz\cap[\Omjpr,\Omje)\cup\{0\}$, 
we have $i_\al=i_\eta$, $\al_0=\al[k-1]^{i_\eta}=\thti(\De+\eta_0)$, and for all such $\be$ \[\al_0[\be]=\thti(\De+\eta_0[\be])\in\Tquok.\]
\end{enumerate}
\end{lem}
{\bf Proof.} We first consider the 
\\[2mm]
{\bf Special case:} $\al=\om^\om$, that is, $i=0$, $\De=0$, and $\eta=\om$. Then we have $i_\eta=0$, $\eta_0=0$,
$\al[k-1]=\thtnod(k-1)=\om^k$, and $\underline{\al[k-1]}=\thtnod(k-2)\in\Tcircquok$ is verified easily by induction on $k$.
\\[2mm]
{\bf Regular case:} $\al>\om^\om$. We then have $\al[k-1]=\thti(\De+\eta[k-1])$. Let $(\Omi,\al_1,\ldots,\alm)$ be the $\Omi$-localization of $\al$. 
Since the case $\eta\not\in\Hz$ is easier to handle, we assume that $\eta\in\Hz$. 
According to Definition \ref{socletermdefi} we then have $\alcirc=\thti(\De+\uual)$ where $\uual:=\almmin$ if $\De^\stari<\almmin<\eta$ with $m>1$,
and $\uual:=0$ otherwise. 
\\[2mm]
\noindent{\bf Case 1:} $l=0$. Then due to Lemma \ref{localizationlem} we have $\alpr=\al[k-1]^{i_\eta+1}=\thti(\De+\eta_0+k)$, hence $\ualpr=\thti(\De+\eta_0+k-1)$, where by
assumption $\De+\eta_0+k-1\in\Tquok$. 

Iteration of the operation $\cdot[0]$ applied to $\eta_0+k-1$ yields a finite strictly decreasing sequence $(\eta_r)_{1\le r}$ of 
ordinals from $\eta_1:=\eta_0+k-1$ down to $0$, passing through all initial sums of $\eta_0+k-1$. In the case where $i=0$ and we eventually reach an ordinal $\eta_\rstar$ such that 
$F_0(\De,\eta_\rstar)$ holds, we may stop at $\rstar$ and note that then $\eta_\rstar=\almmin$ with $m>1$ and $\alcirc=\thti(\De+\almmin)$. 
Otherwise we have $\alcirc=\thti(\De)$ and let $\rstar$ be such that $\eta_\rstar=(\eta_0+k-1)[0]^\rstar=0$. Descending via $\cdot[0]$ from $\ualpr=\thti(\De+\eta_0+k-1)$ passes through the
ordinals $\thti(\De+\eta_r)$ for $r=1,\ldots,\rstar$, where by assumption $\thti(\De+\eta_\rstar)=\alcirc\in\Tquok$. 
For $r\in[1,\rstar)$ we then either have $\eta_r\in\Lim$ such that $F_i(\De,\eta_r)$ does not hold, in which case we have $\thti(\De+\eta_r)[0]=\thti(\De+\eta_{r+1})$ with $\eta_{r+1}=\eta_r[0]$,
or $\eta_r=\eta_{r+1}+1$, in which case we define $\xi_r$ to be the support term of $\thti(\De+\eta_r)$, that is, $\xi_r:=\thti(\De+\eta_{r+1})$. 
For $r\in[1,\rstar)$ we thus have    
\[\thti(\De+\eta_r)[0]=\left\{\begin{array}{cl}
\thti(\De+\eta_{r+1}) & \mbox{ if }\eta_r\in\Lim, \mbox{ otherwise:}\\[2mm]
\xi_r & \mbox{ if }\De=0\\[2mm]
\thti(\De[\xi_r]) & \mbox{ if }\chiomie(\De)=1\\[2mm]
\thti(\De[0]+\xi_r) & \mbox{ otherwise.}
\end{array}\right.\]
For $r=\rstar-1$ we therefore have $\thti(\De+\eta_r)[0]=\alcirc\in\Tquok$ if either $\eta_r\in\Lim$ or $\De=0$.
For $r\in[1,\rstar)$ such that $\eta_r=\eta_{r+1}+1$ and $\De>0$ the support term $\xi_r$ will be preserved when descending further via $\cdot[0]$, thereby reducing the fixed point level $\De$, 
until we reach fixed point level $0$ and hence $\xi_r$. 

Consider the situation where $\eta_r\not\in\Lim$ and $\chiomie(\De)=1$. 
If $\De=\Depr+\Omie$ for some $\Depr$, so that $\De[0]=\Depr$, we obtain $\thti(\De+\eta_r)[0]=\thti(\Depr+\xi_r)$ where $F_i(\Depr,\xi_r)$ holds. 
Otherwise, using Lemma \ref{zeroiterationlem}, $\thti(\De[\xi_r])[0]=\thti(\De[\xi_r[0]]+\xi_r)$ 
since $\De[\xi_r][0]^\stari<\De[\xi_r]^\stari=\xi_r$ with $\xi_r$ the (non-trivial) predecessor of $\thti(\De[\xi_r])$ in
its $\Omi$-localization, so that $\underline{\thti(\De[\xi_r])}=\xi_r$. According to Lemma \ref{zeroiterationlem}, iterating the application of $\cdot[0]$ eventually yields the term 
$\thti(\De[0]+\xi_r)$ in case of $i>0$, and $\thti(\De[1]+\xi_r)$ in case of $i=0$. Terms passed intermediately are of the form $\thti(\De[\xi_r[0]^h]+\xi_r)$.
Note that in case of $i=0$ we have $\domf(\De)=1$ and $\De[1]\in\Tquok$ according to Lemma \ref{Tquoksubstlem}, and descending further via $\cdot[0]$ follows the same pattern with 
constant support term $\xi_r$, until the descent from $\De$ via $\cdot[0]$ and $\cdot[1]$ (occurring occasionally if $i=0$) eventually leads to fixed-point level $0$, followed by $\xi_r$ itself. 

Note that the ordinals above $\alcirc$ in $\zeroseq(\ualpr)$ are contained in the interval $(\almmin,\al)$.
We have to verify that none of these ordinals can be of the form $\ga[l+1]$ where $\ga\in\Tcirc\cap\Lim$ and $l\ge k-2$. 
Since each of these ordinals $\ualpr[0]^h$ with $h$ such that $\ualpr[0]^h>\alcirc$ is a $\thti$-term the argument of which is an element of $\Tquok$, if $\ga[l+1]$ were to match such ordinal, 
we would have $l=k-2$ and $\ga$ would be of the form $\thti(\Ga+\rho)$ where $\chiomie(\Ga)=1$ and either $\rho\not\in\Lim$ or $F_i(\Ga,\rho)$ would hold. 
We have seen earlier, cf.\ Lemmas \ref{fixpointinversionlem} and \ref{zeroiterationinversionlem}, that this is impossible if a support term $\xi_r$ is used in the form $\ualpr[0]^h=\thti(\Depr+\xi_r)$ 
or $\ualpr[0]^h=\thti(\Depr[\xi_r])$ where $\Depr<\De$, as such support terms have fixed-point level $\ge\De$, while $\ga[k-2]$ would have a fixed-point level strictly below $\De$.
We are therefore left with $\ualpr[0]^h$ of the form $\thti(\De+\eta_r)$ where $r\in[1,\rstar)$. Suppose first that $\eta_r=\eta_{r+1}+1$. Then $\thti(\De+\eta_r)$ clearly does not match
$\ga[k-1]$. Otherwise, $\eta_r\in\Lim$. This would imply that $\Ga=\De+\Omie$ and $\eta_r=\ga[k-2]$. Thus $\Ga>\De$, hence $\ga>\al$ as $\ga\in(\almmin,\al)$ would be impossible 
by Proposition \ref{localipic}. Since $\al$ is $k$-residual, there is a $\de\in\Tcircquok$ and $j\le i_\de$ such that $\al=\de[k-1]^j$. We then would have
\[\almmin<\ualpr[0]^h=\ga[k-1]<\al=\de[k-1]^j<\ldots<\de,\]
and Bachmann property would entail $\de<\ga$, which yields the contradiction $\de\not\in\Tquok$ invoking Lemma \ref{intervalcollapsinglem}. 
\\[2mm]
\noindent{\bf Case 2:} $l>0$. Then due to Lemma \ref{localizationlem} we have $i_\al=i_\eta$, $\al_0=\al[k-1]^{i_\eta}=\thti(\De+\eta_0)$, and $\al_0[\be]=\thti(\De+\eta_0[\be])$,
where by assumption $\De+\eta_0[\be]\in\Tquok$. Note that $\eta_0>\almmin$ if $\uual>0$, and that \[\almmin\le\alcirc<\al_0.\]

We have to verify that $\al_0[\be]$ is an element of $\Tquok$. By assumption we have $\chiomje(\eta_0)=1$, which implies that $j<\lv(\sumend(\eta_0))\le i$. 
If $\be\in\Lim$, then $\eta_0[\be]$ is a limit ordinal such that 
$F_i(\De,\eta_0[\be])$ does not hold, and according to Lemma \ref{zeroiterationlem} we obtain \[\al_0[\be][0]^h=\thti(\De+\eta_0[\be[0]^h])=\al_0[\be[0]^h]\] as long as 
$\be[0]^{h^\prime}\in\Lim$ for all $h^\prime<h$, reaching eventually for maximal such $\hstar$ 
\[\alpr:=\al_0[\be][0]^\hstar=\left\{\begin{array}{cl} 
\thti(\De+\eta_0[1])=\al_0[1] & \mbox{ if } j=0\\[2mm]
\thti(\De+\eta_0[0])=\al_0[0] & \mbox{ if } j>0.
\end{array}\right.\] 

Due to Corollary \ref{Bachmanncor} ordinals of the form $\al_0[\be]$ where $\be\in\Hz$ can not be of the form $\ga[l+1]$ for any $\ga\in\Tcirc\cap\Lim$, noting that $\al_0\not\in\Tcirc$. 
Clearly, $\alcirc\in\zeroseq(\alpr)$, and we are left with checking that none of the elements of $\zeroseq(\alpr)$ in the interval $(\alcirc,\alpr]$ can be of the form $\ga[l+1]$ where 
$\ga\in\Tcirc\cap\Lim$ and $l\ge k-2$.

If $\eta_0$ is a successor-multiple of $\Omje$, that is, $\eta_0=_\NF\etapr+\Omje$ for some  $\etapr$, then we have $\eta_0[0]=\etapr$ and $\eta_0[1]=\etapr+1$ (relevant only if $j=0$).
And if $\eta_0$ is a limit-multiple of $\Omje$, then $\eta_0[0]$ is a non-zero multiple of $\Omje$ (least example: $\eta_0=\Om_{j+1}^2$) and, in case of $j=0$, $\eta_0[1]$ is a 
limit-multiple of $\Om_1$ where $\domf(\eta_0[1])\le l$ by Lemma \ref{zeroiterationlem} and, according to Lemma \ref{localizationlem}, $F_i(\De,\eta_0[1])$ does not hold, so that
$\alpr[0]=\thti(\De+\eta_0[1][0])$.
\\[2mm]
{\bf Subcase 2.1:} $j>0$. We are done if $\alpr=\alcirc$, so we may assume that $\alpr>\alcirc$. We then have \[\alpr[0]^h=\thti(\De+\eta_0[0]^{h+1})\] for $h<\hstar$ where $\hstar$ is (minmal)
such that $\eta_0[0]^\hstar=\uual$. Since $\alpr[0]^{\hstar-1}=\alcirc\in\Tquok$, for $h<\hstar-1$ we have $\eta_0[0]^h\in\Lim$ such that $F_i(\De,\eta_0[0]^h)$ does not hold, and 
assuming that $\alpr[0]^h=\ga[l+1]$ we argue as in Case 1,
as the assumption would imply that $\ga=\thti(\Ga+\rho)$ where either $\rho\not\in\Lim$ or $F_i(\Ga,\rho)$ holds, and where $\Ga=\De+\Omie$, $l=k-2$, and $\eta_0[0]^h=\ga[k-2]$. 
\\[2mm]
{\bf Subcase 2.2:} $j=0$. We have already seen that $\alpr=\al_0[1]=\thti(\De+\eta_0[1])$ can not be of the form $\ga[l+1]$. 
\\[2mm]
{\bf 2.2.1:} $\eta_0[1]=\etapr+1$ for some $\etapr\ge0$ such that $\Omje\mid\etapr$ and $\eta_0=_\NF\etapr+\Omje$.  We then have \[\alcirc\le\ualpr=\al_0[0]=\thti(\De+\etapr)\] and
\[\alpr[0]=\left\{\begin{array}{cl}
\ualpr & \mbox{ if }\De=0\\[2mm]
\thti(\De[\ualpr]) & \mbox{ if }\chiomie(\De)=1\\[2mm]
\thti(\De[0]+\ualpr) & \mbox{ otherwise.}
\end{array}\right.\]
As we have seen before, terms $\alpr[0]^{h+1}$ above $\ualpr=\al_0[0]$ can not be of the form $\ga[l+1]$. And for $\alpr[0]^{h+1}\le\al_0[0]$ we argue as in Subcase 2.1. 
\\[2mm]
{\bf 2.2.2:} $\eta_0[1]$ is a limit-multiple of $\Omje$.
Then let $\hstar$ be (minimal) such that $\eta_0[1][0]^\hstar=\uual$. Then, for $h\le\hstar$ we have \[\alpr[0]^h=\thti(\De+\eta_0[1][0]^h).\]
Here we argue as in Subcase 2.1 with $\eta_0[1]$ in place of $\eta_0[0]$.
\qed

\begin{lem}\label{extendedTquoksubstlem}
Let $\al=\thti(\De+\eta)$, where $\De>0$ with $\chiomie(\De)=0$ and either $\eta\not\in\Lim$ or $F_i(\De,\eta)$ holds. Suppose that $\al$ is $k$-residual such that 
$\almmin, \ual, \De\in\Tquok$ where $(\Omi,\al_1,\ldots,\alm)$ is the $\Omi$-localization of $\al$. 
Assume further that for $\De_0=\De[k-1]^{i_\De}$ we have $\De_0[\be]\in\Tquok$ for $\be\in\Tquok\cap\Hz\cap[\Omjpr,\Omje)\cup\{0\}$ where $j$ satisfies $\domf(\De_0)=j+1$.
Using the abbreviation $\alpr:=\al[k-1]^{i_\De}$, we have $\ualpr\in\Tquok$ and distinguish between the following cases for $j$: 
\begin{enumerate}
\item If $j<i$, we have $i_\al=i_\De$, $\alpr=\al_0$, and for all $\be$ as above \[\al_0[\be]\in\Tquok.\]
\item If $j=i$ (and hence $i_\De>0$), we have $\alpr=\thti(\De_0+\uual)$ where $\uual\in\{0,\almmin,\ual\}$, 
\[\alpr[0],\ldots,\alpr[k-2]\in\Tquok,\]
and therefore $\alpr\in\Tfcircquok$.
\end{enumerate}
\end{lem}
{\bf Proof.} 
We begin with calculating the ordinals $\alpr:=\al[k-1]^{i_\De}$ and $\ualpr$. 
In case of $j<i$ we have $i_\al=i_\De$, $\alpr=\al_0$, which is equal to $\al$ in the special case $i_\De=0$, and also calculate $\al_0[\be]$.  
For these calculations we distinguish between situations $\ual>0$ and $\ual=0$:

If $\ual>0$, it follows that $\alpr=\thti(\De_0+\ual)$, unless $\alpr=\al$ in case of $i_\De=0$. We have $\ualpr=\ual$ and, in case of $j<i$, $\al_0[\be]=\thti(\De_0[\be]+\ual)$. 
Note that if $\alpr<\al$, we have $\De_0^\stari<\ual$ by parts \ref{starcontrolclaim} and \ref{stardomclaim} of Lemma \ref{bracketsmainlem}.

If $\ual=0$ (and hence $\eta=0$), note that according to Lemma \ref{localizationlem} $\almmin$ is the immediate predecessor of $\alpr$ in its $\Omi$-localization.
We prove by induction on $i_\De$ that either $\alpr=\thti(\De_0)$ or $\alpr=\thti(\De_0+\al_{m-1})$, where the latter applies if $m>1$ and
\begin{equation}\label{uualeq}
\De[k-1]^l[0]^\stari<(\De[k-1]^l)^\stari=\almmin
\end{equation}
for a least $l< i_\De$, which, if yes, implies that $l>0$ as $\ual=0$, and $\De[k-1]^{l-1}[0]^\stari=\almmin\le\De^\stari$ according to Corollary \ref{deepsupporttermcor}. 
Similarly, in case of $j<i$, either $\al_0[\be]=\thti(\De_0[\be])$ or $\al_0[\be]=\thti(\De_0[\be]+\al_{m-1})$, 
depending on whether $m>1$ and (\ref{uualeq}) holds, this time for some (least) $l\le i_\De$. We also see that $\ualpr\in\{0,\almmin\}$ according to the latter criterion.

Indeed, if $i_\De=0$ and hence $\De_0=\De$ and $j<i$ as $\chiomie(\De)=0$ by assumption, we have $\al_0=\thti(\De_0)=\al$, thus $\ualpr=0$
and $\al_0[\be]=\thti(\De_0[\be])$. 

If $i_\De>0$, we have $\al[k-1]=\thti(\De[k-1])$ as $\ual=0$.
Taking a look specifically at the case $i_\De=1$, we have $\alpr=\al[k-1]=\thti(\De_0)$. If $j=i$ and hence $\chiomie(\De_0)=1$, we have $\ualpr=0$. If $j<i$, we have $\ualpr=0$, unless
$\De_0[0]^\stari<\De_0^\stari=\almmin$ with $m>1$, using Lemma \ref{localizationlem}. 
In this latter case $\ualpr=\almmin=\De[0]^\stari\le\De^\stari$ since $\almmin\le\De[0]^\stari\le\De[k-1]^\stari=\De_0^\stari\le\De^\stari$ according to Corollary \ref{deepsupporttermcor} and 
part \ref{starcontrolclaim} of Lemma \ref{bracketsmainlem}. Thus, $\al_0[\be]=\thti(\De_0[\be]+\ualpr)$ if $j<i$.  

If $i_\De>1$, let $l\le i_\De$ be minimal such that $m>1$ and (\ref{uualeq}) holds, if that exists, and $l:=i_\De+1$ otherwise. If $l\ge i_\De$, we immediately see that 
$\alpr=\thti(\De_0)$, and if $j<i$, we have $\ualpr=\almmin$ if $l=i_\De$, otherwise $\ualpr=0$, and $\al_0[\be]=\thti(\De_0[\be]+\ualpr)$. Suppose from now that $l<i_\De$. Then we have
$m>1$ and \[\De^\stari\ge\De[k-1]^\stari\ge\ldots\ge(\De[k-1]^l)^\stari=\almmin>\De[k-1]^l[0]^\stari\] using part \ref{starcontrolclaim} of Lemma \ref{bracketsmainlem},
so that \[\al[k-1]^\lpr=\left\{\begin{array}{cl}
\thti(\De[k-1]^\lpr) & \mbox{ for }\lpr\le l, \mbox{ and}\\[2mm] 
\thti(\De[k-1]^\lpr+\almmin) & \mbox{ for }\lpr\in(l, i_\De],
\end{array}\right.\]
thus $\alpr=\thti(\De_0+\almmin)$ and $\ualpr=\almmin$. If $j<i$, we have $\alpr=\al_0$ and $\al_0[\be]=\thti(\De_0[\be]+\almmin)$. Due to Corollary \ref{deepsupporttermcor} for $\lpr<l$ we have 
\[\De[k-1]^\lpr[0]^\stari\ge\almmin,\] with equality for $\lpr=l-1$, while clearly $(\De[k-1]^\lpr)^\stari<\almmin$ for $\lpr>l$, using part \ref{stardomclaim} of Lemma \ref{bracketsmainlem}.

So far, we have seen that $\ualpr\in\{0,\almmin,\ual\}$ and hence $\ualpr\in\Tquok$, and proceed with the verification of the two claims of the lemma.
\\[2mm]
{\bf Ad 1:} Suppose that $j<i$. Then $\De_0$ is a limit-multiple of $\Omie$, and we have $\al_0[\be]=\thti(\De_0[\be]+\ualpr)$.

If $\ualpr>0$, then $F_i(\De_0[\be],\ualpr)$ holds, and Lemmas \ref{fixpointinversionlem} and \ref{zeroiterationinversionlem} show that $\al_0[\be][0]^h$ can not be of the form $\ga[l+1]$ where 
$\ga\in\Tcirc\cap\Lim$ and $l\ge k-2$ for all $h$ such that $\al_0[\be][0]^h>\ualpr$, while clearly $\ualpr\in\zeroseq(\al_0[\be])\cap\Tquok$.

We now consider the situation $\ualpr=0$. Clearly, if $\be>0$, $\al_0[\be]=\thti(\De_0[\be])$ is not of the form $\ga[l+1]$ with $\ga, l$ as above, since $\be<\Omi\le\ga[k-2]$ and due to 
Corollary \ref{Bachmanncor}. 
Recall that according to Lemma \ref{zeroiterationlem} $\al_0[\be][0]=\thti(\De_0[\be[0]])$ for $\be\in\Lim$.
We are left with showing that terms of $\zeroseq(\al_0[1])$ in case of $j=0$, or of $\zeroseq(\al_0[0])$ in case of $j>0$ can not be of the form $\ga[l+1]$. As seen before, of these terms
only those that are strictly above $\almmin$ and that are not using $\almmin$ as explicit support term, still need to be checked. Inspecting with the aid of Lemma \ref{supportbyzeroiterationlem}
shows that such $\al_0[\be][0]^h$ are of the form $\thti(\Deprnod[0]^{r+1})$ where $\Deprnod:=\De_0[1]$ if $j=0$ and $\Deprnod:=\De_0$ if $j>0$, which in case of $m>1$ satisfy 
$(\Deprnod[0]^{r+1})^\stari\ge\almmin$, and which therefore are elements of the interval $(\almmin,\al)$. 
Assuming that $\al_0[\be][0]^h=\ga[k-1]$ implies that $\Deprnod[0]^r[0]=\Ga[\ga[k-2]]$. This entails $\Deprnod[0]^r<\Ga$, which is seen as follows.
Since $\Omie\mid\Deprnod[0]^{r+1}$, it follows that $\Ga$ is a limit-multiple of $\Omie$. We either have $\Deprnod[0]^r=\Deprnod[0]^{r+1}+\Omie<\Ga$, or $\Deprnod[0]^r$ is a limit-multiple
of $\Omie$, and due to Corollary \ref{Bachmanncor} we must have $\Deprnod[0]^r<\Ga$ as assuming the contrary would imply $\ga[k-2]=0$, which is not the case.
Considering the chain of inequalities 
\[\Deprnod[0]^{r+1}<\ldots<\Deprnod[0]=\left\{\begin{array}{clcl} \De_0[0] & < & \De_0[\Omj] & \mbox{ if }j>0\\ \De_0[1][0] & < & \De_0[1] & \mbox{ if } j=0\end{array}\right\}
<\ldots<\De_0[\be]<\De_0=\De[k-1]^{i_\De}<\ldots<\De,\] 
Bachmann property implies that we therefore even have $\De<\Ga$ and hence $\al<\ga$ as $\ga\in(\almmin,\al)$ is impossible by Proposition \ref{localipic}. 
As $\al$ is $k$-residual, there exists $\de\in\Tcircquok$ such that $\al=\de[k-1]^s$ for some $s\le i_\de$, so that, moreover, $\de<\ga$, again due to Bachmann property.
This, however, means that $\al_0[\be][0]^h=\ga[k-1]<\de<\ga$, which according to Lemma \ref{intervalcollapsinglem} leads to the contradiction $\de\not\in\Tquok$. 
This completes the proof of $\al_0[\be]\in\Tquok$.
\\[2mm]
{\bf Ad 2:} Suppose that $j=i$. According to the assumptions we then have $i_\De>0$ and $\chiomie(\De_0)=1$, as well as $\alpr=\thti(\De_0+\uual)$ where $\uual\in\{0,\almmin,\ual\}$.
$\alpr[q]\in\Tquok$ for $q=0,\ldots,k-2$ is shown by induction on $q$. 
\\[2mm]
{\bf Base case:} $q=0$. We have $\alpr[0]=\thti(\De_0[\ualpr])$ where $\ualpr\in\{0,\almmin,\ual\}$ as shown above.
\\[2mm]
{\bf Subcase 1:} $\ualpr>0$. Note that if $\De_0=\De_0[0]+\Omie$, we have $\alpr[0]=\thti(\De_0[0]+\ualpr)$, and if $\De_0$ is a limit-multiple of $\Omie$, according to Lemma \ref{zeroiterationlem}
we have $\alpr[0]^2=\thti(\De_0[\ualpr[0]]+\ualpr)$, etc. As shown in Lemmas \ref{fixpointinversionlem} and \ref{zeroiterationinversionlem}, the elements of $\zeroseq(\alpr[0])$ above $\ualpr$ can not
be of a form $\ga[l+1]$ where $\ga\in\Tcirc\cap\Lim$ and $l\ge k-2$. We therefore have $\alpr[0]\in\Tquok$.
\\[2mm]
{\bf Subcase 2:} $\ualpr=0$. Here we have $\alpr[0]=\thti(\De_0[0])$ and argue as in the proof of claim 1 for $\al_0[0]$.
\\[2mm]
{\bf Inductive step:} Suppose that $\alpr[0],\ldots,\alpr[q]\in\Tquok$ where $q<k-2$. Then $\alpr[q+1]\in\Tquok$ is shown exactly as in the inductive step of the proof of 
Lemma \ref{Tquoksimplefixplem}.
\qed

As an example motivating the following main lemma and illustrating a comparison to the original Goodstein principle in Lemma \ref{alternativelem}, we consider the
special case of $\al=\om^\om=\thtnod(\thtnod(0))$, which is the situation $i=0$, $\De=0$, and $\eta=\om$ in the lemma below. $\al$ is equal to $\thtnod(\thte(0))[1]=\epsn[1]$ 
and therefore not an element of $\Ttwo$ but a $2$-residual of $\epsn\in\Ttwo$.  
We have $i_\eta=0$, $\eta_0=0$, $\al[k-1]=\thtnod(k-1)=\om^k$, and
\[\al[k-1]^2=\om^{k-1}\cdot k,\quad\ldots,\quad\al[k-1]^{k+1}=\om^{k-1}\cdot(k-1)+\ldots+\om\cdot(k-1)+k.\] 
We thus obtain $i_\al=k$ and \[\al_0=\om^{k-1}\cdot(k-1)+\ldots+\om\cdot(k-1)\in\Tcircquok.\]
Note that in the lemma below this situation is covered by Case 1, followed by iterations of Case 2.1.

The main lemma is formulated for additive principal ordinals $\al$, however, note that once the lemma is established, its claim holds for all $\al\in\Tcircquok$
since we are descending along fundamental sequences that target the last summand of an approximated limit ordinal.
\begin{lem}[Main Lemma]\label{mainlem} For $k\ge2$ and $\al\in\Tcircquokcls$, where $\al=\thti(\De+\eta)$,
\begin{enumerate}
\item if $i>0$, for any $\be\in\Tquok\cap\Hz\cap[\Omjpr,\Omje)\cup\{0\}$ where $j+1=\domf(\al_0)$, we have \[\al_0[\be]\in\Tquok,\]
\item if $i=0$, we have \[\al_0\in\Tcircquok.\]
\end{enumerate}
\end{lem}
{\bf Proof.} The proof is by main induction on $i_\al$ and side induction on the build-up of $\al$. 
According to Lemma \ref{supporttermlem} we have $\al\in\Tcirc$ and $\alcirc,\ual\in\Tquok$.
Clearly, we have $\lv(\sumend(\al_0))=i$ and hence $\domf(\al_0)\le i$ so that defining $j$ in case of $i>0$ by $\domf(\al_0)=j+1$, we have $j<i$. 
Since in the cases where $\al_0=0$ or where $\al_0$ is a successor-multiple of $\Omje$ we are done immediately, we may assume that $\al_0\in\Lim$ or, in case of 
$i>0$, that $\al_0$ is a limit-multiple of $\Omje$.
We distinguish between two main cases:
\\[2mm]
\noindent{\bf Case 1:} $\eta\in\Lim$ and $F_i(\De,\eta)$ does not hold. 
Then we have $\domf(\al)=\domf(\eta)$, and $\al[k-1]=\thti(\De+\eta[k-1])$, unless $i=0$, $\De=0$, and $\eta=\om$, in which case we have $\al[k-1]=\thtnod(k-1)$.
$\al$ is $k$-residual, and either $\al\in\Tcircquok$ or $\al=\ga[k-1]\in\Tlcircquok$ for some $\ga\in\Tfcircquok$, hence in particular $\eta\in\Tcircquok$ (where in the latter situation
$\eta=\ga[k-2]\in\Hz^{\ge\Omi}$) and $\De+\eta\in\Tquok$. 
The i.h.\ applies to $\eta$ (more precisely: to $\sumend(\eta)$) as we see that $i_\eta<i_\al$, so that the assumptions of Lemma \ref{Tquoklimlem} are fulfilled: if $\lv(\sumend(\eta))=0$, we obtain 
$\alpr:=\al[k-1]^{i_\eta+1}\in\Tscircquok$ and apply the i.h.\ to $\alpr$. Otherwise we also have $i>0$, and part 2 of Lemma \ref{Tquoklimlem} yields $\al_0[\be]\in\Tquok$ immediately.
\\[2mm]
\noindent{\bf Case 2:} $\eta\not\in\Lim$ or $F_i(\De,\eta)$ holds.
\\[2mm]
\noindent{\bf Subcase 2.1:} $\De=0$. If $\eta=0$ we have $\al=\Omi$ so that $i=0$ and $\al=\om$, hence $\al[k-1]=k$ and $\al_0=0$.
In case of $F_i(0,\eta)$ we have $\ual=\eta$, while $\ual=\thti(\etapr)$ if $\eta=\etapr+1$ for some $\etapr$. We have $\al[k-1]=\ual\cdot k$, where $\ual\in\Tquok$ by Lemma \ref{supporttermlem},
and if $\ual\in\Tcircquok$, it follows that $\ual\cdot(k-1)\in\Tcircquok$, and we may apply the i.h.\ to $\ual$ as $i_{\ual}<i_\al$.
Otherwise $\al_0=\al[k-1]$, and Lemma \ref{Tquoksubstlem} yields $\al_0[\be]=\ual\cdot(k-1)+\ual[\be]\in\Tquok$. 
\\[2mm]
\noindent{\bf Subcase 2.2:} $\chiomie(\De)=1$. We then have either $\al\in\Tcircquok\cup\Tscircquok$, or $\al\in\Tfcircquok$.
For $\al\in\Tcircquok\cup\Tscircquok$, according to Lemma \ref{Tquoksubstlem} we have 
$\De[\be]\in\Tquok$ for any $\be\in\Tquok\cap\Hz\cap[\Omipr,\Omie)\cup\{0\}$, and Lemma \ref{Tquoksimplefixplem} yields $\al[0],\ldots,\al[k-2]\in\Tquok$.
Therefore, $\al\in\Tfcircquok$ and hence $\al[0],\ldots,\al[k-2]\in\Tquok$  in any case, $(\al[k-1])^\circ\in\Tquok$ by Lemma \ref{supporttermlem}, and if $\al[k-1]\in\Tcirc$ 
and hence $\al[k-1]\in\Tlcircquok$, we may apply the i.h.\ to $\al[k-1]$. 

Otherwise we have $\al_0=\al[k-1]$ with $\domf(\al_0)=j+1$ for some $j<i$ and need to show that $\al_0[\be]\in\Tquok$ for $\be\in\Tquok\cap\Hz\cap[\Omjpr,\Omje)\cup\{0\}$.
By Lemma \ref{zeroiterationlem} we have $\domf(\al_0)=\domf(\De[\al[k-2]])=\domf(\al[k-2])$ and $\al_0[\be]=\thti(\De[\al[k-2][\be]])$, where using Lemma \ref{Tquoksubstlem} 
$\al[k-2][\be]\in\Tquok$ and hence $\De[\al[k-2][\be]]\in\Tquok$. We consider cases regarding $\De$.
\begin{enumerate}
\item $\De=\Depr+\Omie$. Then we have $\al[k-1]=\thti(\Depr+\al[k-2])$ where $\al[k-2]\in\Hz\cap[\Omi,\Omie)$ and $F_i(\Depr,\al[k-2])$ does not hold. 
Here part 2 of Lemma \ref{Tquoklimlem} applies to $\al_0=\thti(\Depr+\al[k-2])$, i.e.\ with $\al[k-2]$ taking the role of $\eta=\eta_0$.
\item $\De$, hence also $\De[\al[k-2]]$, is a limit-multiple of $\Omie$. Here Lemma \ref{extendedTquoksubstlem} part 1 applies to $\al_0=\thti(\De[\al[k-2]])$.
\end{enumerate}
\noindent{\bf Subcase 2.3:} $\De$ is a limit-multiple of $\Omie$ such that $\chiomie(\De)=0$. Then $\domf(\De)=\domf(\al)=0$, and for $\al$ one of the following cases applies:
\begin{itemize}
\item $\al\in\Tcircquok\cup\Tscircquok$, or 
\item $\al=\ga[k-1]\in\Tcirc$ where $\ga=\thti(\Ga+\rho)\in\Tfcircquok$, so that $\De=\Ga[\ga[k-2]]\in\Tquok$, $\eta=0$, and $\ual=0$ as shown in the proof of Lemma \ref{supporttermlem}. 
\end{itemize}
Let $(\Omi,\al_1,\ldots,\al_m)$ be the $\Omi$-localization of $\al$. 
In any case we have $\almmin, \alcirc, \ual,\De\in\Tquok$, and the i.h.\ applies to $\De$ as we see that $i_\De\le i_\al$ with $\De$ a proper subterm of $\al$. 
Thus $\De_0=\De[k-1]^{i_\De}$ with $\chiomje(\De_0)=1$ for some $j\le i$ and $\De_0[\be]\in\Tquok$ for $\be\in\Tquok\cap\Hz\cap[\Omjpr,\Omje)\cup\{0\}$. 
The assumptions of Lemma \ref{extendedTquoksubstlem} are therefore fulfilled:
if $j<i$, then $\al_0[\be]\in\Tquok$ follows from part 1 of Lemma \ref{extendedTquoksubstlem},
and if $j=i$, the i.h.\ applies to $\alpr:=\al[k-1]^{i_\De}$, as $\alpr\in\Tfcircquok$ according to part 2 of Lemma \ref{extendedTquoksubstlem}.  
\qed

\begin{cor}\label{maincor} For $k\ge2$ the predecessor function $\Pk$ restricted to $\Tk$ is the true predecessor operation on $\Tk$.
\end{cor}
{\bf Proof.} Note that for $\al\in\Tk\cap\Lim$ we have $\Pk\al=\al_0+k-1$. 
\qed

\begin{theo}\label{maintheo} For $k\ge2$ the restriction of $\Ck$ to the set $\Tk$ is an order isomorphism with image $\N$. We thus have \[\Tk=\TcapOmquok.\]
\end{theo}
{\bf Proof.} Using Lemmas \ref{intervalcollapsinglem}  and \ref{mainlem} we argue as in the proof of Corollary \ref{esccollapscor} and apply Lemma \ref{easyinclusionlem} and
Corollary \ref{maincor}.  
 \qed

\section{Goodstein process}\label{Goodsteinprocesssec}

\begin{defi}\label{basechangedefi} Let $N\in\N$ and $k\in[2,\om)$. According to Theorem \ref{maintheo}  there exists a unique $\al\in\Tk$ such that $N=\Ck(\al)$,
the base-$k$ representation of $N$. 
For any $l\in(k,\om)$ the change of base from $k$ to $l$ for $N$ is defined by \[N[k\mapsto l]:=\Cl(\al),\]
where $\al$ is the base-$l$ representation of $N[k\mapsto l]$ since $\Tk\subseteq\Tl$. 
The change of base $k$ to $\om$ for $N$ is defined by \[N[k\mapsto\om]:=\al.\]

As outlined in Remark \ref{uniformityrmk} we proceed in the same way in the context of a notation system $(\tau,\cdot[\cdot])$, relying on order isomorphisms
of $\tau/_{\Ck}$ (in place of $\Tk$) and $\N$ as in Theorem \ref{maxquotisomtheo} with $\tau$ in place of $\T\cap\Om_1$, cf.\ \cite{W}.
\end{defi}

Now that Theorems \ref{maxquotisomtheo} and \ref{maintheo} are established, Cichon's trick applies as layed out in the introduction, and we obtain the reformulation 
of the generalized Goodstein principle (\ref{reformulation}). Since the function defined in (\ref{hkdefieq}) is expressed in terms of Hardy hierarchy as in (\ref{hkHardyeq}), 
if the generalized Goodstein principle were provable in $\pioneonecanod$, we would obtain a proof of the totality of $H_\tau$ for $\tau:=\T\cap\Om_1$ in 
$\pioneonecanod$ using the canonical ``master'' fundamental sequence given by $\theta_n:=\thtnod(\ldots(\thtn(0))\ldots)$, which is completely contained in each $\Tk$, 
so that $H_\tau(k)=H_{\theta_k}(k)$, despite the well-known fact that $H_\tau$ is not provably total in $\pioneonecanod$. In case of an abstractly given notation system
$(\tau,\cdot[\cdot])$ as specified in Remark \ref{uniformityrmk} we proceed similarly using maximality quotients $\tau/_{\Ck}$. We thus obtain

\begin{cor} The Goodstein process defined according to base transformation as in Definition \ref{basechangedefi} always terminates, which however, when based on 
ordinals below $\T\cap\Om_1$, exceeds the proof-strength of the theory $\pioneonecanod$. 

A uniform result is obtained for any theory $\mathcal{T}$ that allows for a proper ordinal analysis with 
proof-theoretic ordinal $\tau$, for which we have a system $(\tau,\cdot[\cdot])$ of fundamental sequences with Bachmann property and a 
suitable ``master'' sequence as specified in Remark \ref{uniformityrmk}, and according to which $H_\tau$ is not provably total in $\mathcal{T}$, cf.\ \cite{W}.
\end{cor}

\section{Conclusion and outlook}

The concept of quotients allows for a most
direct approach that establishes an immediate connection between complexity of ordinal notations and fast growing number theoretic functions. A natural next 
step would be to consider quotients of restricted iteration of a notation system from \cite{BS88} of strength $\pioneonetrnod$, in which uniqueness of terms is achieved 
via normal form conditions, and for which Bachmann property was shown in \cite{We95}.

Recent progress related to Goodstein principles was achieved by \cite{AFDWW20} and \cite{AWW21}, elaborating on Arai's idea of normal form representations for 
natural numbers using the Ackermann function (and variants) and also building on earlier ideas by Weiermann in \cite{We17}. That approach was further extended in 
\cite{FDW24b} and \cite{FDW25} to a wider variety of Goodstein principles for which termination and independence can be proved uniformly via the notion 
of \textit{base-change maximality} of normal form representations of natural numbers. 

We see the present work as a basis for the deeper analysis of hierarchies of fast-growing functions and a greater variety of Goodstein principles in the spirit 
of \cite{AFDWW20,AWW21,FDW24b,FDW25}.

\section*{Acknowledgement}
I would like to thank Professor Andreas Weiermann for numerous insightful and motivating discussions on ordinal notation systems and Goodstein principles,
and for the suggestion to verify the characterization in terms of maximality in Corollary \ref{maximalitycor}.

\end{document}